\documentclass{agtart_a}
\pdfoutput=1

\usepackage{graphicx}
\usepackage{dsfont}
\usepackage{color}

\title[The plastikstufe -- a generalization of the overtwisted disk]{The plastikstufe -- a generalization of the overtwisted disk to higher dimensions}

\author{Klaus Niederkruger}
\givenname{Klaus}
\surname{Niederkruger}
\address{Departement de Mathematiques\\
Universite Libre de Bruxelles, CP 218\\\newline
Boulevard du Triomphe\\
B-1050 Bruxelles\\Belgium}
\email{kniederk@ulb.ac.be}
\urladdr{}

\volumenumber{6}
\issuenumber{}
\publicationyear{2006}
\papernumber{86}
\startpage{2473}
\endpage{2508}

\doi{}
\MR{}
\Zbl{}


\received{11 September 2006}
\revised{17 November 2006}
\accepted{19 November 2006}
\published{15 December 2006}
\publishedonline{15 December 2006}
\proposed{}
\seconded{}
\corresponding{}
\editor{Liz}
\version{}

\arxivreference{math.SG/0607610}

\AtBeginDocument{\let\bar\wbar\let\tilde\wtilde\let\hat\what}
\AtBeginDocument{\renewcommand{\phi}{\varphi}
  \renewcommand{\epsilon}{\varepsilon}
  \renewcommand{\theta}{\vartheta}
  \renewcommand{\S}{{\mathbb{S}}}}

\def\SetFigFont#1#2#3#4#5{\small} 
\def\adjustlabel<#1,#2>#3{\smash{\rlap{\kern #1 \raise #2\hbox{#3}}}}

\newcommand{\N}{{\mathbb{N}}}

\newcommand{\1}{{\mathds{1}}}

\newcommand{\T}{{\mathbb{T}}}
\newcommand{\Disk}[1]{{\mathbb{D}^{#1}}}
\newcommand{\Reeb}{{X_\mathrm{Reeb}}}
\newcommand{\lie}[1]{{\mathcal{L}_{#1}}}

\newcommand{\plastikstufe}[1]{{\mathcal{PS}({#1})}}
\newcommand{\norm}[1]{{\lVert #1\rVert}}
\newcommand{\abs}[1]{{\left\lvert #1\right\rvert}}
\newcommand{\overtwisted}{{\mathbb{D}_\mathrm{OT}}}
\newcommand{\lcan}{{\lambda_{\mathrm{can}}}}
\newcommand{\POS}{{\mathbf{q}}}
\newcommand{\MOM}{{\mathbf{p}}}

\DeclareMathOperator{\End}{End}
\DeclareMathOperator{\evaluationmap}{ev}
\DeclareMathOperator{\id}{id}
\DeclareMathOperator{\ImaginaryPart}{Im}
\DeclareMathOperator{\linkingzahl}{lk}

\DeclareMathOperator{\RealPart}{Re}
\DeclareMathOperator{\SO}{SO}

\makeatletter
\def\cnewtheorem#1[#2]#3{\newtheorem{#1}{#3}
\expandafter\let\csname c@#1\endcsname\c@theorem}
\makeatother

\theoremstyle{plain}
\newtheorem{theorem}{Theorem}
\cnewtheorem{propo}[theorem]{Proposition}
\cnewtheorem{coro}[theorem]{Corollary}
\newtheorem{question}{Question}

\theoremstyle{remark}
\newtheorem{remark}{Remark}
\newtheorem{example}{Example}

\theoremstyle{definition}
\newtheorem*{defi}{Definition}

\begin{document}

\begin{asciiabstract}
  In this article, we give a first prototype-definition of
  overtwistedness in higher dimensions. According to this definition,
  a contact manifold is called ``overtwisted'' if it contains a
  ``plastikstufe'', a submanifold foliated by the contact structure
  in a certain way.  In three dimensions the definition of the
  plastikstufe is identical to the one of the overtwisted disk.  The
  main justification for this definition lies in the fact that the
  existence of a plastikstufe implies that the contact manifold does
  not have a (semipositive) symplectic filling.
\end{asciiabstract}

\begin{abstract}
  In this article, we give a first prototype-definition of
  overtwistedness in higher dimensions. According to this definition,
  a contact manifold is called \textit{overtwisted} if it contains a
  \textit{plastikstufe}, a submanifold foliated by the contact structure
  in a certain way.  In three dimensions the definition of the
  plastikstufe is identical to the one of the overtwisted disk.  The
  main justification for this definition lies in the fact that the
  existence of a plastikstufe implies that the contact manifold does
  not have a (semipositive) symplectic filling.
\end{abstract}

\maketitle
The situation of contact topology can be roughly stated like this: the
$3$--dimensional contact manifolds can be understood very adequately
by topological methods, a far-reaching classification has been
achieved and relations to many other fields have been established.
In contrast, the world map of higher-dimensional contact geometry
consists almost entirely of white spots.  A powerful method for
constructing such manifolds is \textit{contact surgery}, the most
promising technique developed so far to distinguish different contact
structures is \textit{contact homology} and with Giroux's
\textit{open book decomposition}, it is hoped that some classification
results could be obtained.

The first structural distinction found for contact $3$--manifolds was
the notion of overtwistedness.  It turned out that such manifolds
firstly do not allow an (even weak) symplectic filling
by Eliashberg \cite{Eliashberg_HoloDiscs} and Gromov \cite{Gromov_Kurven}, and secondly can be
classified in a very satisfactory way as in Eliashberg \cite{Eliashberg_Overtwisted}.

In higher dimensions, surprisingly, no analogous criterion has yet
been found.  Giroux has proposed a definition based on his open
book decomposition, which in three dimensions is completely equivalent
to the standard one.  In contrast, our definition is based on the
existence of a \textit{plastikstufe}, a direct generalization of the
overtwisted disk.  In Gromov's famous paper on holomorphic
curves \cite{Gromov_Kurven}, a sketchy description of something, which
possibly could be a plastikstufe, is given.  The generalization of
overtwistedness described in this article was found independently by
Yuri Chekanov.  Interestingly, his (unpublished) proof of
\fullref{hauptsatz} uses very different methods.

The definition of the plastikstufe given in this paper is certainly
only a preliminary version, meant as a prototype leading to a
criterion for nonfillability in higher dimensions.  Our
definition implies the following theorem.

\begin{theorem}\label{hauptsatz}
  Let $(M,\alpha)$ be a contact manifold containing an embedded
  plastikstufe.  Then $M$ does not have any semipositive symplectic
  filling.  If $\dim M \le 5$, then $M$ does not have any symplectic
  filling at all.
\end{theorem}

The proof of this statement will be given in \fullref{sec: beweis
  hauptsatz}.

\begin{remark}
  A $2n$--dimensional symplectic manifold $(M,\omega)$ is called
  \textit{semipositive} if every $A\in\pi_2(M)$ with $\omega(A) > 0$
  and $c_1(A)\ge 3 - n$ has nonnegative Chern number.  Note that
  every symplectic $4$-- or $6$--manifold is semipositive.
\end{remark}

There are several shortcomings of the definition of overtwisted given
here, the most important being that no example of a \textit{closed}
contact manifold containing an embedded plastikstufe has been found so
far.  It is relatively easy though to construct open manifolds
containing an embedded plastikstufe.  As observed by Chekanov,
the plastikstufe can be used to detect exotic contact structures on
$\R^{2n-1}$ (see \fullref{sec: exotische kontakstruktur}).

From a practical viewpoint the definition of the plastikstufe is also
rather cumbersome, because it is less topological than the overtwisted
disk in dimension~$3$ (\fullref{bemerkungen ueber frobenius}).

\begin{remark}
  At the time of the final revision of this article, Francisco
    Presas Mata announced a method which allows one to construct
  closed contact manifolds that contain embedded plastikstufes
  \cite{PresasExamplesPlastikstufes}.
\end{remark}

\paragraph{Acknowledgments}

The work on this article was initiated at the Universit\"at zu
K\"oln and finished at the Universit\'e Libre de Bruxelles.  My
research was funded in Germany by the university and in Belgium by
the Fonds National de la Recherche Scientifique (FNRS).

The definition for the plastikstufe given in this article was found
after innumerable discussions with Kai Zehmisch.  Later, it was
Fr\'ed\'eric Bourgeois who guided me in converting the intuitive
picture I had in mind into sound and hard mathematics.   I spoke 
with Otto van Koert during the whole project almost on a
daily basis.  These conversations helped me find solutions to many
problems in this article.  Without these three colleagues, this
article wouldn't have been possible. Furthermore, Francisco
Presas Mata helped me solve the last steps in the proof.  I also
profited from discussions with Yuri Chekanov, Hansj\"org
Geiges, Ferit \"Ozt\"urk, Dietmar Salamon, Felix
Schlenk and Sava\c{s} Yaz{\i}c{\i}.

\setcounter{section}{-1}

\section{Preliminaries}
\label{sec: wiederholung von fakten}

The following notions are standard in symplectic topology, but for
completeness we briefly repeat them here.

\begin{defi}
  Let $(W,\omega)$ be a symplectic manifold.  A \textit{Liouville
    vector field} $X_L$ is a vector field on $W$, whose flow makes the
  symplectic form expand exponentially. This property can be
  formulated equivalently as
  \begin{equation*}
    \lie{X_L}\omega = \omega.
  \end{equation*}
  A \textit{(convex) symplectic filling} $(W,\omega)$ of a contact
  manifold $(M,\alpha)$ is a compact symplectic manifold with boundary
  $\partial W = M$, such that there exists a Liouville vector field
  $X_L$ in a neighborhood of $\partial W$ that points outwards, and
  such that
  \begin{equation*}
    \alpha = \left.\bigl(\iota_{X_L}\omega\bigr)\right|_{TM}.
  \end{equation*}
\end{defi}

\begin{remark}
  We can define a function $h\co U\to(-\infty,0]$ on a neighborhood
  $U\subset W$ of $M$, by considering the time $t_p$ it takes a point
  $p\in U$ to flow along $X_L$ to $M$, and then setting $h(p) :=
  -t_p$.  Define $\tilde\alpha := \iota_{X_L}\omega$.  By taking on
  each level set $M_t:=h^{-1}(t)$ (sufficiently close to $M$) the Reeb
  field of the contact form $\alpha_t := \left. \tilde\alpha
  \right|_{TM_t}$, we obtain a smooth vector field, which we denote by
  $\Reeb$.
\end{remark}

In the context of this article we will use the term ``compatible
almost complex structure'' in the following sense.

\begin{defi}
  Let $(W,\omega)$ be a symplectic filling of a contact manifold
  $(M,\alpha)$.  A compatible almost complex structure $J$ is a smooth
  section of the endomorphism bundle $\End(TW)$ such that $J^2
  = -\1$, that is compatible with $\omega$ in the usual sense, which
  means that for all $X,Y\in T_pW$, the following equation holds
 \begin{gather*}
   \omega(JX,JY) = \omega(X,Y),\\
   \tag*{\hbox{and}}
   g(X,Y) := \omega(JX,Y)
 \end{gather*}
 defines a Riemannian metric.  Additionally, we require $J$ to satisfy
 close to the boundary $M = \partial W$ the following properties: for
 the two vector fields $X_L$ and $\Reeb$ introduced above, $J$ is
 defined as
 $$
   JX_L = \Reeb \quad \text{and}\quad J\Reeb = - X_L,
 $$
 and $J$ leaves the subbundle $\xi_t = \ker \alpha_t \le TM_t$
 invariant.
\end{defi}

\begin{propo}
  Let $u\co V\cap \mathbb{H} \to W$ be a $J$--holomorphic map
  ($V\subset\C$ is an open set, and $\mathbb{H}\subset\C$ is the upper
  halfplane).  The function $h\circ u\co V\cap \mathbb{H} \to \R$ is
  subharmonic.
\end{propo}
\begin{proof}
  A short computation shows that
$$ \eqalignbot{
    0 & \le u^*\omega = u^*d\iota_{X_L}\omega = u^*d\tilde \alpha =
    u^* d \bigl(- dh\circ J\bigr) = -u^* dd^c h \cr
    & = -dd^c(h\circ u) = \Bigl(\frac{\partial^2 h\circ u}{\partial
      x^2}+ \frac{\partial^2 h\circ u}{\partial y^2} \Bigr)\,dx\wedge
    dy.\cr 
} \proved $$
\end{proof}

\begin{coro}\label{kurven transvers zu rand} By the strong maximum
  principle and the boundary point lemma (eg\ Gilbarg and Trudinger
  \cite{GilbargTrudinger}), any $J$--holomorphic curve
  $u\co (\Sigma,\partial\Sigma) \to (W,\partial W)$ is either constant
  or it intersects $M = \partial W$ only at $\partial\Sigma$, and this
  intersection is transverse.
\end{coro}

Finally, we will denote the moduli space of $J$--holomorphic disks
lying in $U$ with boundary in $\tilde U$ and with one marked point
$z_0\in\partial\Disk{2}$ by the symbol $\mathcal{M}(U,\tilde U, z_0)$.

\section{Definition of the plastikstufe}\label{sec: definition der plastikstufe}

Before giving the definition of the plastikstufe, we first introduce
some preliminary definitions.

\begin{defi}
  A \textit{maximally foliated submanifold} $L$ in a
  $(2n-1)$--dimensional contact manifold $(M,\alpha)$ is a submanifold
  of dimension $n$ on which $\ker\left.\alpha\right|_{TL}$ defines a
  (possibly singular) foliation.
\end{defi}

\begin{remark}\label{bemerkungen ueber maximal geblaettert} The term
  ``maximally'' in the definition above means that $L$ is not
  contained in some higher dimensional submanifold also foliated by
  $\alpha$.  The condition on the dimension is imposed by the fact
  that the leaves of the foliation are locally Legendrian
  submanifolds.
\end{remark}
\begin{remark}\label{bemerkungen ueber frobenius}
  Frobenius' Theorem implies that $L$ is foliated by $\left.\alpha\right|_{TL}$
  if and only if $\left.(\alpha\wedge d\alpha)\right|_{TL} \equiv 0$.
\end{remark}
\begin{remark}
  The main reason why these submanifolds are interesting in the
  setting of this paper, is that if $(M,\alpha)$ is the convex
  boundary of a symplectic manifold $(W,\omega)$, and $W$ is given an
  almost complex structure $J$ compatible on $M$ with $\xi =
  \ker\alpha$, then a maximally foliated submanifold $L\hookrightarrow
  M$ will be (at the nonsingular points of the foliation) a totally
  real submanifold in $W$ such that the Fredholm theory of
  $J$--holomorphic curves can be applied.
\end{remark}

\begin{defi}\label{definition elliptische singularitaet} An
  \textit{elliptic singular set} $S$ inside a maximally foliated
  submanifold $L$ is a closed codimension~$2$ submanifold inside $L$
  whose neighborhood is diffeomorphic to $\Disk{2}\times S
  \hookrightarrow L$ with coordinates $(x,y;s)$ such that
  $\left.\alpha\right|_{TL}$ is represented by $x\,dy - y\,dx$ on this
  neighborhood (see \fullref{bild: elliptische singularitaet}).
\end{defi}

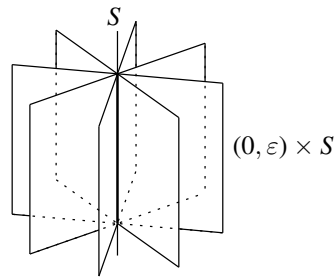
\begin{figure}[ht!]
  \begin{center}
    \begin{picture}(0,0)%
      \includegraphics{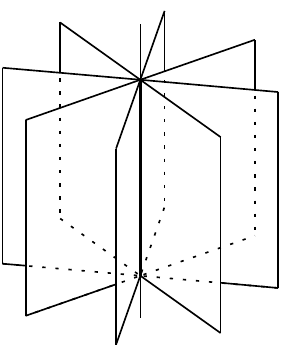}%
    \end{picture}%
    \setlength{\unitlength}{4144sp}%
    \begingroup\makeatletter\ifx\SetFigFont\undefined%
    \gdef\SetFigFont#1#2#3#4#5{%
      \reset@font\fontsize{#1}{#2pt}%
      \fontfamily{#3}\fontseries{#4}\fontshape{#5}%
      \selectfont}%
    \fi\endgroup%
    \begin{picture}(1319,1585)(-1,-1025)
      \put(1318,-281){\makebox(0,0)[lb]{\smash{{\SetFigFont{7}{8.4}{\familydefault}{\mddefault}{\updefault}{\color[rgb]{0,0,0}$(0,\epsilon)\times S$}%
            }}}}
      \put(567,491){\makebox(0,0)[lb]{\smash{{\SetFigFont{7}{8.4}{\familydefault}{\mddefault}{\updefault}{\color[rgb]{0,0,0}$S$}%
            }}}}
    \end{picture}%
    \caption{The foliation around a set of elliptic singularities
      consists of a circle of stripes $(0,\epsilon)\times S$ forming
      a circle of rays touching in the singular set.}\label{bild:
      elliptische singularitaet}
  \end{center}
\end{figure}

\begin{example}
  It is very easy to find examples of maximally foliated submanifolds
  with elliptic singularities (at least locally).  As
  \fullref{bild: elliptische singularitaet} already suggests, there
  is a similarity between an elliptic singularity and the binding of
  an open book.
  
  By a result of Giroux \cite{Giroux}, every contact manifold
  $(M,\alpha)$ has a compatible open book decomposition $(P,\theta)$
  with binding $(B,\alpha_B)$, where $B=\partial P$ and $\alpha_B =
  \left.\alpha\right|_{TB}$.  The normal form of the neighborhood of
  $(B,\alpha_B)$ can then be chosen to be
  \begin{align*}
    \Bigl(B\times \Disk{2}, \alpha_B + \frac{1}{2}\,(x\,dy -
    y\,dx)\Bigr),
  \end{align*}
  where $(x,y)$ are the coordinates on the $2$--disk.  If one chooses
  any Legendrian submanifold $S$ inside $B$ (which exist in abundance
  by Ekholm, Etnyre and Sullivan \cite{EkholmEtnyre}), then the set $S\times\Disk{2} \hookrightarrow
  B\times\Disk{2}$ is a maximally foliated submanifold with elliptic
  singular set $S$.
  
  An interesting problem would then be to try extend this submanifold
  into $M-B$.
\end{example}

The following definition is fundamental in $3$--dimensional contact
topology.

\begin{defi}
  Let $(M,\alpha)$ be a $3$--dimensional contact manifold. An embedded
  $2$--disk
  $$
  \iota\co \Disk{2} \hookrightarrow M
  $$
  is called an \textit{overtwisted disk} $\overtwisted$, if there is
  only one point on the disk where the foliation given by
  $\iota^*\alpha$ is singular, and if the boundary of $\overtwisted$
  is the only closed leaf of this foliation (see \fullref{bild:
    blaetterung auf ueberdrehter scheibe}).
\end{defi}

\begin{figure}[ht!]
  \begin{center}
    \begin{picture}(0,0)%
      \includegraphics{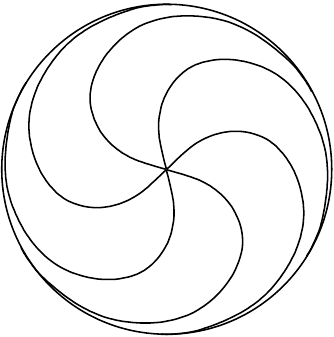}%
    \end{picture}%
    \setlength{\unitlength}{3947sp}%
    \begingroup\makeatletter\ifx\SetFigFont\undefined%
    \gdef\SetFigFont#1#2#3#4#5{%
      \reset@font\fontsize{#1}{#2pt}%
      \fontfamily{#3}\fontseries{#4}\fontshape{#5}%
      \selectfont}%
    \fi\endgroup%
    \begin{picture}(1600,1608)(2393,-2550)
    \end{picture}%
    \caption{Foliation induced by $\iota^*\alpha$ on the overtwisted
      disk}\label{bild: blaetterung auf ueberdrehter scheibe}
  \end{center}
\end{figure}

Now we will give a conceivable generalization to higher dimensions:
let $(M,\alpha)$ be a $(2n-1)$--dimensional contact manifold, and let
$S$ be a closed $(n-2)$--dimensional manifold.
\begin{defi}
  A \textit{plastikstufe $\plastikstufe{S}$ with singular set $S$} in
  $M$ is an embedding of the $n$--dimensional manifold
  $$
  \iota\co \Disk{2}\times S \hookrightarrow M
  $$
  that is maximally foliated by the $1$--form $\beta := \iota^*\alpha$.
  The boundary $\partial\plastikstufe{S}$ of the plastikstufe should
  be the only closed leaf, and there should be an elliptic singular
  set at $\{0\}\times S$.  The rest of the plastikstufe should be
  foliated by an $\S^1$--family of stripes, each one diffeomorphic to
  $(0,1)\times S$, which are spanned between the singular set on one
  end and approach $\partial\plastikstufe{S}$ on the other side
  asymptotically.
\end{defi}

\begin{remark}\label{bemerkungen ueber plastikstufe} An overtwisted
  disk $\overtwisted$ is equal to a $2$--dimensional plastikstufe
  $\plastikstufe{\{p\}}$.
\end{remark}

\begin{remark}
  As mentioned in \fullref{bemerkungen ueber frobenius},
  $\ker\beta$ defines a foliation on $\Disk{2}\times S$, if and only
  if $\beta\wedge d\beta \equiv 0$.  This means that the definition
  above requires $\iota$ to satisfy a partial differential equation,
  in contrast to $3$--dimensional contact topology, where the foliation
  condition is trivially satisfied.
\end{remark}

\begin{remark}
  The boundary of the plastikstufe is a Legendrian submanifold
  $\partial\plastikstufe{S} \cong \S^1\times S$ of~$M$.
\end{remark}

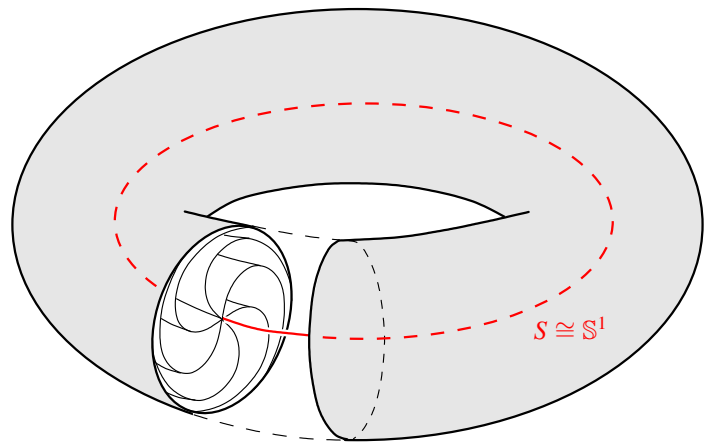
\begin{figure}[ht!]
  \begin{center}
    \begin{picture}(0,0)%
      \includegraphics[scale=.9]{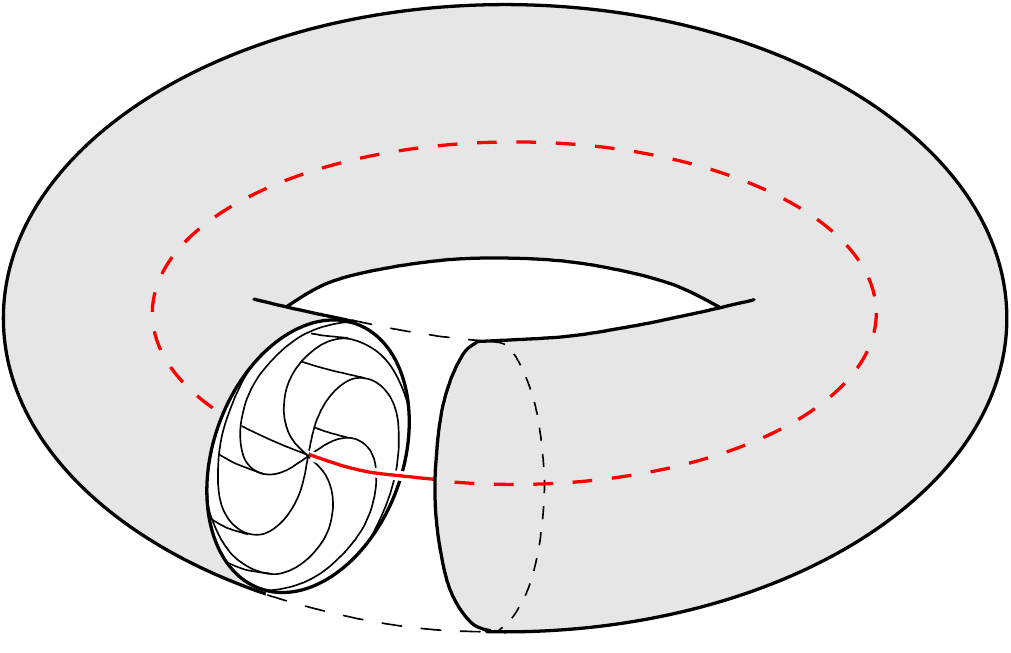}%
    \end{picture}%
    \setlength{\unitlength}{3730sp}%
    \begingroup\makeatletter\ifx\SetFigFont\undefined%
    \gdef\SetFigFont#1#2#3#4#5{%
      \reset@font\fontsize{#1}{#2pt}%
      \fontfamily{#3}\fontseries{#4}\fontshape{#5}%
      \selectfont}%
    \fi\endgroup%
    \begin{picture}(4618,2978)(1843,-8212)
      \put(5311,-7441){\makebox(0,0)[lb]{\smash{{\SetFigFont{14}{16.8}{\familydefault}{\mddefault}{\updefault}{\color[rgb]{1,0,0}$S\cong\S^1$}%
            }}}}
    \end{picture}%
    \caption{The plastikstufe in a $5$--dimensional contact
      manifold}\label{bild: plastikstufe in dimension 5}
  \end{center}
\end{figure}

\begin{defi}
  A $(2n-1)$--dimensional contact manifold $(M,\alpha)$ will be called
  \textit{overtwisted} if it contains an $n$--dimensional plastikstufe.
\end{defi}

\begin{example}
  It is easy to construct for any closed manifold $S$ an
  (unfortunately only) open contact manifold that contains the
  plastikstufe $\plastikstufe{S}$: let $(N^3,\alpha_0)$ be an
  overtwisted $3$--manifold with overtwisted disk
  $\iota_0\co\overtwisted\hookrightarrow N$. Let $M$ be the
  $(2n-1)$--dimensional manifold $M := N \times T^*S$ with the
  $1$--form
  \begin{equation*}
    \alpha = \alpha_0 + \lcan,
  \end{equation*}
  where $\lcan = -\MOM\cdot d\POS$ is the canonical $1$--form of the
  cotangent bundle $T^*S$.
  
  The $1$--form $\alpha$ is a contact form, and it is easy to check
  that the embedding
  \begin{equation*}
    \iota\co\plastikstufe{S}\hookrightarrow M,\, (p,s) \mapsto
    (\iota_0(p),\sigma_0(s)),
  \end{equation*}
  where $\sigma_0\co S\hookrightarrow T^*S$ is the zero-section, really
  defines an $n$--dimensional plastikstufe.
\end{example}

\section[Sketch of the proof of Theorem 1]{Sketch of the proof of \fullref{hauptsatz}}\label{sec: beweis hauptsatz}

The proof follows the line of the original proof in dimension~$3$.  We
study a certain moduli space $\mathcal{M}$ of holomorphic disks, and
there we find a cycle representing the trivial homology class in
$\mathcal{M}$, but at the same time we show that it is mapped by the
evaluation map to a nontrivial element of homology thus producing a
contradiction.

Let $(M,\alpha)$ be a $(2n-1)$--dimensional closed compact manifold
that contains a plastikstufe $\plastikstufe{S}$.  Assume that
$(W,\omega)$ is a semipositive symplectic filling of~$M$.  We will
choose a compatible almost complex structure $J$ on $W$ (in the sense
of the definition given in \fullref{sec: wiederholung von
  fakten}).
The most important properties of $J$--holomorphic curves with respect
to the convex boundary $M$ is that a holomorphic curve $u\co\Sigma \to
W$ cannot be tangent to $M$, and in particular if $p\in u(\Sigma)$ is
a point, where the holomorphic curve intersects $M$, then $T_pu\cap
\xi_p = \{0\}$.  Note that $T_p u = 0$ at a point $p\in M$ implies
that $u\equiv p$.

In the rest of the proof, we consider the $J$--holomorphic disks
$u\co\Disk{2}\to W$, whose boundary $u(\partial\Disk{2})$ lies on the
plastikstufe, and which have a marked point $z_0\in\partial\Disk{2}$
on their boundary.  If $u\not\equiv \mathrm{const}$, then
$u(\partial\Disk{2})$ is transverse to the foliation on
$\plastikstufe{S}$, and $u(\partial\Disk{2})$ cannot touch the
boundary of $\plastikstufe{S}$.  The boundary of every such curve is
linked with the singular set $S$, because otherwise it could not be
transverse to the foliation.  This also implies that the only constant
curves in the compactification of the moduli space lie in $S$.

In \fullref{sec: bishop family}, we study the Bishop family of
disks emanating from the singular set of the plastikstufe: for a small
neighborhood $U$ of $S$, we find a standard model, where we can choose
a carefully prepared almost complex structure $J$ (which can be
extended to a regular $J$ on the whole symplectic filling).  For this
$J$ and for $\tilde U:= U\cap \plastikstufe{S}$, the evaluation map
$\evaluationmap_{z_0}$ is a diffeomorphism between the moduli space
$\mathcal{M}(U,\tilde U,z_0)$ and $\tilde U$.  Stated differently, for
every point $p\in \tilde U$ there is a unique holomorphic disk $u$
lying inside $U$ such that $u(z_0) = p$.  We can explicitly write
down this continuous family
\begin{align*}
  \psi\co & \tilde U \to \mathcal{M}(W,\plastikstufe{S},z_0),
\end{align*}
such that $\evaluationmap_{z_0}\circ \psi = \id_{\tilde U}$.  The
curves $\psi(s)$ for every $s\in S$ are constant disks, and the
boundary of every other disk $\psi(p)$ with $p\notin S$ is linked once
with the singular set $S$.

Next, we show that the model neighborhood $U$ is foliated by compact
codimension~$2$ manifolds, which are $J$--holomorphic.  An intersection
argument then allows us to show that any holomorphic disk $u$ lying
partially in $U$ and whose boundary can be capped off by attaching a
disk lying inside $\plastikstufe{S}$ has to be completely contained in
$U$.  Hence the disk $u$ lies in the image of the map $\psi$.

With the results obtained in \fullref{sec: bishop family}, one can
see that the $N$--dimensional moduli space of holomorphic disks ($N =
\dim\plastikstufe{S}$) described above only has a single end touching
the singular set $S$, and stays otherwise at a finite distance from
$S$.

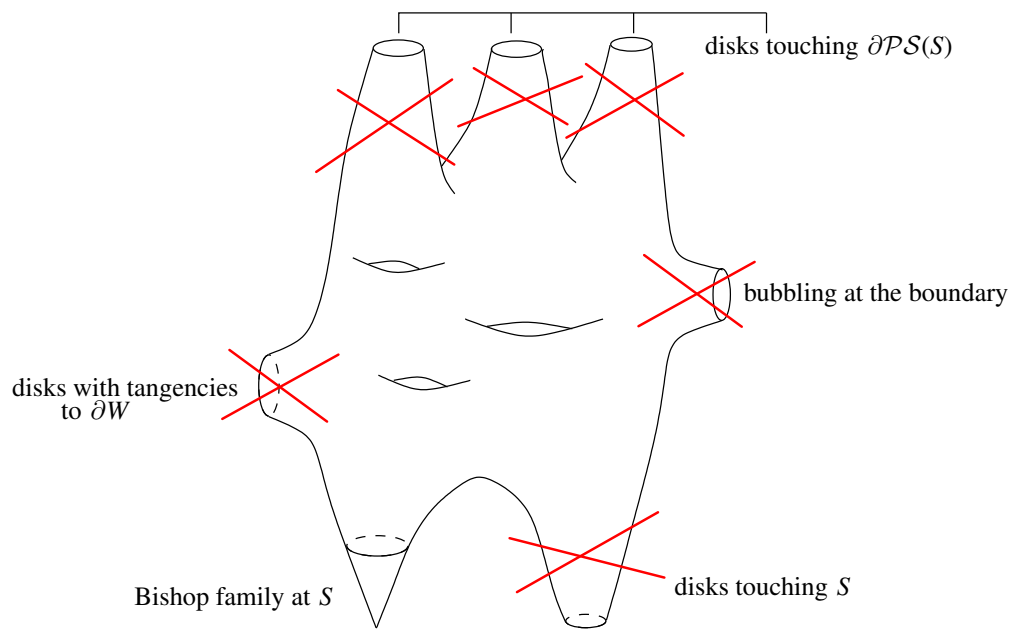
\begin{figure}[ht!]
  \begin{center}
    \begin{picture}(0,0)%
      \includegraphics{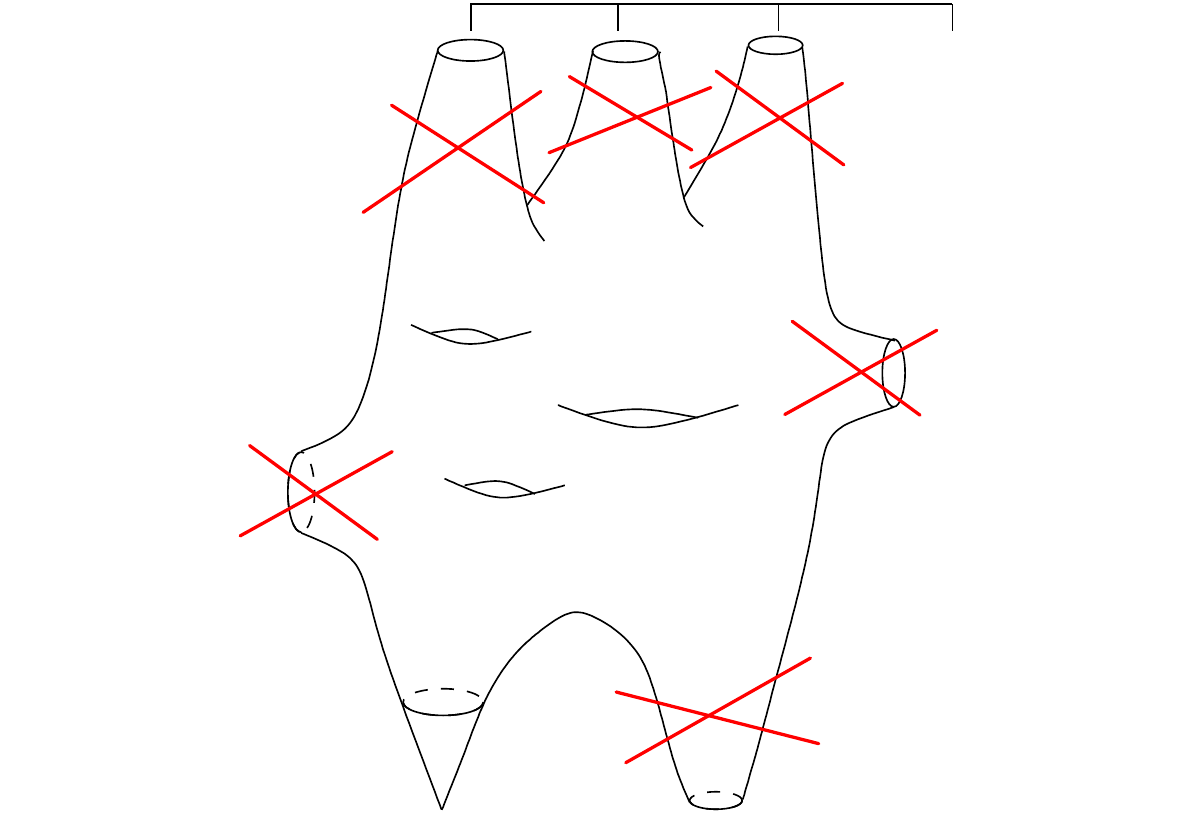}%
    \end{picture}%
    \setlength{\unitlength}{4144sp}%
    \begingroup\makeatletter\ifx\SetFigFont\undefined%
    \gdef\SetFigFont#1#2#3#4#5{%
      \reset@font\fontsize{#1}{#2pt}%
      \fontfamily{#3}\fontseries{#4}\fontshape{#5}%
      \selectfont}%
    \fi\endgroup%
    \begin{picture}(5465,3708)(541,-7807)
      \put(4528,-4355){\makebox(0,0)[lb]{\smash{{\SetFigFont{8}{9.6}{\familydefault}{\mddefault}{\updefault}{\color[rgb]{0,0,0}disks touching $\partial\plastikstufe{S}$}%
            }}}}
      \put(4345,-7594){\makebox(0,0)[lb]{\smash{{\SetFigFont{8}{9.6}{\familydefault}{\mddefault}{\updefault}{\color[rgb]{0,0,0}disks touching $S$}%
            }}}}
      \put(1351,-7655){\makebox(0,0)[lb]{\smash{{\SetFigFont{8}{9.6}{\familydefault}{\mddefault}{\updefault}{\adjustlabel<-15pt,0pt>{\color[rgb]{0,0,0}Bishop family at $S$}}%
            }}}}
      \put(4760,-5834){\makebox(0,0)[lb]{\smash{{\SetFigFont{8}{9.6}{\familydefault}{\mddefault}{\updefault}{\color[rgb]{0,0,0}bubbling at the boundary}%
            }}}}
      \put(541,-6406){\makebox(0,0)[lb]{\smash{{\SetFigFont{8}{9.6}{\familydefault}{\mddefault}{\updefault}{\adjustlabel<-10pt,0pt>{\color[rgb]{0,0,0}disks with tangencies}}%
            }}}}
      \put(676,-6541){\makebox(0,0)[lb]{\smash{{\SetFigFont{8}{9.6}{\familydefault}{\mddefault}{\updefault}{\color[rgb]{0,0,0}to $\partial W$}%
            }}}}
    \end{picture}%
    \caption{The possible boundary components of the moduli space
      $\mathcal{M}(W,\plastikstufe{S},z_0)$ are given by disks
      touching the singular set $S$ or the boundary of the
      plastikstufe $\partial\plastikstufe{S}$, disks with tangencies
      to the boundary of the symplectic manifold $\partial W$, or
      curves having bubbles at their boundary. Different arguments
      show that with exception of the Bishop family all of these cases
      can be excluded.}\label{bild: moduli raum mit rand}
  \end{center}
\end{figure}

The compactness result of Gromov states that the moduli space
of holomorphic disks with bounded energy is always compact, provided
we allow bubbling.  More precisely, in a compact symplectic manifold
every moduli space of simple disks with uniformly bounded energy,
whose boundary sits on a compact totally real submanifold is a smooth
manifold that can be compactified by including bubbled curves.  In
\fullref{sec: no bubbling}, we first show that there is a uniform
energy bound for all curves in $\mathcal{M}(W,\plastikstufe{S},z_0)$.
In our situation, $\plastikstufe{S} - S$ is not a compact totally real
manifold, but since only the Bishop family comes close to $S$, the
compactness property still holds in a suitable sense.  There are two
possible types of bubbles that can occur: either a sphere can bubble
off at the interior of the holomorphic disks or a new holomorphic disk
can form at the boundary of the disks.  The foliation of the
plastikstufe $\plastikstufe{S}$ imposes a constraint on holomorphic
curves, which forbids the second type of bubbling.

\begin{figure}[ht!]
  \begin{center}
    \begin{picture}(0,0)%
      \includegraphics{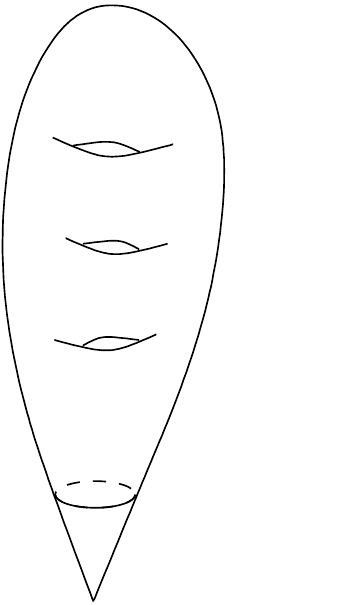}%
    \end{picture}%
    \setlength{\unitlength}{4144sp}%
    \begingroup\makeatletter\ifx\SetFigFont\undefined%
    \gdef\SetFigFont#1#2#3#4#5{%
      \reset@font\fontsize{#1}{#2pt}%
      \fontfamily{#3}\fontseries{#4}\fontshape{#5}%
      \selectfont}%
    \fi\endgroup%
    \begin{picture}(1645,2750)(2131,-7804)
      \put(2746,-7666){\makebox(0,0)[lb]{\smash{{\SetFigFont{8}{9.6}{\familydefault}{\mddefault}{\updefault}{\color[rgb]{0,0,0}Bishop family at $S$}%
            }}}}
      \put(2161,-7351){\makebox(0,0)[lb]{\smash{{\SetFigFont{8}{9.6}{\familydefault}{\mddefault}{\updefault}{\color[rgb]{0,0,0}$C_\epsilon$}%
            }}}}
    \end{picture}%
    \caption{The moduli space $\mathcal{M}(W,\plastikstufe{S},z_0)$
      only has one boundary component, and hence it follows that
      $[C_\epsilon]$ represent the trivial class in
      $H_{N}\bigl(\mathcal{M}(W,\plastikstufe{S},z_0)\bigr)$ for $N =
      \dim\mathcal{M} -1$.}\label{bild: moduli raum ohne rand}
  \end{center}
\end{figure}
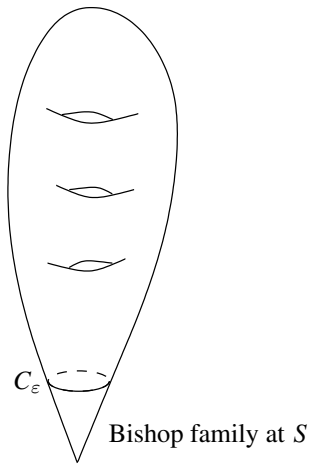

The treatment of interior bubbles is more technical, so we will first
assume that no bubbling at all can happen to illustrate more easily
the geometrical idea of the proof.  The moduli space
$\mathcal{M}(W,\plastikstufe{S},z_0)$ is a smooth manifold, which
could have several ``boundary components'' (see \fullref{bild:
  moduli raum mit rand}).  By the arguments given so far (and by the
assumption that no spheres can bubble off), the only boundary of the
moduli space corresponds to the Bishop family $\psi$ (see
\fullref{bild: moduli raum ohne rand}).

With the evaluation map at $z_0\in\partial\Disk{2}$, it is easy to see
that
\begin{align*}
  C_\epsilon &:= \Bigl\{\psi(p) \Bigm| p \in \plastikstufe{S},\,
  d(p,S) = \epsilon\Bigr\},
\end{align*}
for $\epsilon >0$ sufficiently small, maps to the generator of
$H_{N-1}(\plastikstufe{S}-S)$ with $N=\dim\plastikstufe{S}$.  Hence it
follows that $C_\epsilon$ has to represents a nontrivial cycle in the
moduli space (see \fullref{bild: moduli raum ohne rand}).  But this
leads to a contradiction, because $C_\epsilon$ is the only boundary
component of the moduli space (if we remove all $C_r$ with
$r<\epsilon$), and hence it represents a boundary in homology.

The conclusion is that the Bishop family cannot be the only boundary
component of the moduli space.  The only assumption made to exclude
other components was that $W$ is a convex filling of the contact
manifold $M$.  This assumption must be false, and hence $M$ is not
fillable.

The general situation, ie\ the one where interior bubbles can occur,
can be treated like this: trace a path $\gamma(I)$ (with $I=[0,1]$) on
$\plastikstufe{S}$ running from the singular set $S$ to the boundary
$\partial\plastikstufe{S}$.  If the evaluation map
$\smash{\evaluationmap_{z_0}\co\mathcal{M}(W,\plastikstufe{S},z_0)\to
\plastikstufe{S}}$ is transverse to $\gamma$, then
$\smash{\evaluationmap_{z_0}^{-1}\bigl(\gamma(I)\bigr)}$ will be a smooth
$1$--dimensional submanifold of $\mathcal{M}(W,
\plastikstufe{S},z_0)$.  Define a projection $\pi$ by
\begin{align*}
  \pi\co & \plastikstufe{S} - S \to \partial\plastikstufe{S} \\
  & (z,s) \mapsto \Bigl(\frac{z}{\abs{z}},s\Bigr).
\end{align*}
The preimage of a point $(e^{i\phi_0},s_0)$ under $\pi$ is a radial
line in the plastikstufe joining $S$ and $\partial \plastikstufe{S}$.
By Sard's~Theorem, the set of singular values of the map $\pi\circ
\evaluationmap_{z_0}$ has measure~$0$, and because that map is
surjective (consider the restriction to the Bishop family $\psi$), we
find plenty of radial curves $\gamma$ such that
$\evaluationmap_{z_0}^{-1}\bigl(\gamma(I)\bigr)$ will be a
$1$--dimensional submanifold of $\mathcal{M}(W,\plastikstufe{S},z_0)$.

In \fullref{sec: no bubbling}, we show that the evaluation map
$\evaluationmap_{z_0}$ is a pseudocycle (in the sense of
McDuff and Salamon \cite{McDuffSalamonJHolo}), which means that the image of the
evaluation map on the complement $\wwbar{\mathcal{M}} - \mathcal{M}$ (ie\ on the
bubbled curves) lies in the closed image of finitely many
codimension~$2$ manifolds.  The projection of this set under $\pi$
still has at least codimension~$1$, and hence, almost every point in
$\partial\plastikstufe{S}$ is a regular value of $\pi\circ
\evaluationmap_{z_0}$ and not contained in the image of the bubbled
curves.

For such a value $(e^{i\phi_0},s_0)$, the corresponding preimage in
$\mathcal{M}(W,\plastikstufe{S},z_0)$ is a collection of compact
$1$--dimensional submanifolds that can only have boundary points, at
the boundary of the moduli space.  Since it avoids the bubble curves,
it can actually only have boundary points on $C_\epsilon$.  But
because the Bishop family $\psi$ is a diffeomorphism, the
$1$--dimensional submanifold only touches $C_\epsilon$ exactly once.
This is a contradiction, because it means there is one component of
the preimage that is a compact $1$--dimensional manifold with only one
boundary component.

\section{A Bishop family around elliptic singularities}\label{sec:
  bishop family}

The aim of this section is to show that there exists a unique family
of holomorphic curves growing out of each component $S$ of the
elliptic singular set of a maximally foliated submanifold.  I'm
indebted to Fr\'ed\'eric Bourgeois for orienting me through the
theory of Cauchy--Riemann operators.

The main idea consists in arranging a certain almost complex structure
around~$S$, where it is possible to do all required computations
explicitly.

\subsection{A local model for a neighborhood of the singular set}

Let $(M,\alpha)$ denote always in this section a contact
$(2n-1)$--manifold, which has a symplectic filling $(W,\omega)$.
Assume further that $M$ contains a maximally foliated submanifold $F$,
with an elliptic singular set, and let $S$ be one component of this
set.  To simplify the calculation, we want to find some standard form
for the neighborhood of the singular set $S$.

\begin{propo}\label{standard kontaktform bei singularitaet} There is a
  small neighborhood of the singular set $S\subset F$ in $M$ that is
  contactomorphic to a neighborhood of $\{0\}\times S$ in the contact
  manifold
  \begin{align*}
    \Bigl( \Disk{3}\times T^*S, dz + \frac{1}{2}\,(x\,dy-y\,dx) +
    \lcan\Bigr),
  \end{align*}
  where $(x,y,z)$ are the coordinates on the $3$--ball, and $\lcan =
  -\MOM\cdot d\POS = -\sum_j p_j\,dq^j$ is the canonical $1$--form on
  $T^*S$.  In this neighborhood, the maximally foliated submanifold
  $F$ lies in the set $\bigl\{(x,y,0;\sigma_0(s))\bigr\}$, where
  $\sigma_0$ denotes the zero-section in $T^*S$.
\end{propo}
\begin{proof}
  By our definition of an elliptic singular set, we find a
  neighborhood around $S$ inside $F$ that can be written as
  $\{(w;s)\in\C\times S\}$, and the restriction of the contact form is
  equal to
  \begin{align*}
    \left.\alpha\right|_{TF} &= \frac{1}{2}\,\bigl(x\,dy -
    y\,dx\bigr),
  \end{align*}
  with $w=x+iy$.

  Choose now a $d\alpha$--compatible complex structure $J$ on
  $\xi=\ker\alpha$.  Note that the tangent space of $S_{w_0} :=
  \{w_0\}\times S$ for any fixed $w_0\in\C$ lies in the contact
  structure $\xi$, and secondly that $J\cdot TS_{w_0}$ is transverse
  to $F$, because if there was a nonzero vector $X\in TS_{w_0}$ such
  that $JX\in TF$, then
  \begin{align*}
    0 &\ne d\alpha (X,JX) = dx\wedge dy(X,JX) = 0,
  \end{align*}
  which is a contradiction.  Hence it follows that
  \begin{equation*}
    J\cdot (TS_{w_0}) \cap TF = \{0\}.
  \end{equation*}
  Similarly, since the Reeb field $\Reeb$ is transverse to the contact
  structure $\xi$, it follows in particular that $\Reeb$ is transverse
  both to $J\cdot(TS_{w_0})$, and (at least close to the singular set
  $S$) to the maximally foliated submanifold $F$, because $TF\le\xi$
  on~$S$.

  Choose now a metric $g$ on $M$ such that $J\cdot (TS_{w_0}) \perp
  F$, $\Reeb \perp F$, and $\norm{\Reeb} = 1$.  This gives an
  identification for the normal bundle
  \begin{align*}
    \nu F \cong \langle\Reeb\rangle \oplus J\cdot (TS_{w_0}),
  \end{align*}
  which can be combined with the map
  \begin{align*}
    \langle\Reeb\rangle \oplus J\cdot (TS_{w_0}) &\to \R^3\times T^*S\\
    (x,y,\POS;z\cdot \Reeb + J\dot\POS) & \mapsto (x,y,z;\POS, g(\dot
    \POS,\cdot)).
  \end{align*}
  By using the exponential map as in the proof of the tubular
  neighborhood theorem, one gets a diffeomorphism from a neighborhood
  of $\{0\}\times S$ in $\R^3\times T^*S$ to a neighborhood of the
  singular set $S$ in $M$ such that $E:=\{(x,y,0;\sigma_0(S)\}$ is
  mapped into the submanifold $F$.

  The pullback of the contact form evaluates in this model on $E$ to
  \begin{align*}
    \left.\alpha\right|_E &= dz + \frac{1}{2}\,(x\,dy-y\,dx),
  \end{align*}
  because $\partial_z$ is equal to the Reeb field, and the restriction
  of $\alpha$ to $F$ is equal to the second term.  The $\POS$-- and
  $\MOM$--directions lie at every point of $E$ in the contact
  structure.  The $2$--form $d\alpha$ is written on $E$ as
  \begin{align*}
    \left.d\alpha\right|_E &= dx\wedge dy + d\POS\wedge d\MOM +
    \mathrm{Rest},
  \end{align*}
  where ``$\mathrm{Rest}$'' are terms pairing $dx$ and $dy$ with
  $\MOM$--coordinates.

  In the final step, we use now an improved version of the Moser trick
  (as explained for example in Geiges \cite[Theorem~2.24]{Geiges_Handbook})
  to find a vector field $X_t$ that isotopes the contact form given
  into the desired one $dz + \frac{1}{2}\,(x\,dy-y\,dx) + \lcan$.  Let
  $\alpha_t$, $t\in[0,1]$, be the linear interpolation between both
  $1$--forms.  Assume there is an isotopy $\psi_t$ defined around $S$
  such that $\psi_t^*\alpha_t = \alpha_0$.  The field $X_t$ generating
  this isotopy satisfies the equation
  \begin{equation*}
    \lie{X_t} \alpha_t +\dot \alpha_t = 0.
  \end{equation*}
  By writing $X_t = H_t\,R_t + Y_t$, where $H_t$ is a smooth function,
  $R_t$ is the Reeb vector field of $\alpha_t$, and $Y_t\in\ker
  \alpha_t$, we obtain plugging then $R_t$ into the equation above
  \begin{align*}
    dH_t(R_t) = -\dot \alpha_t(R_t).
  \end{align*}
  The vector field $Y_t$ is completely determined by $H_t$, because
  $Y_t$ satisfies the equations
  \begin{align*}
    \iota_{Y_t}\alpha_t &= 0, \\
    \iota_{Y_t}d\alpha_t &= - dH_t - \dot \alpha_t,
  \end{align*}
  hence it suffices to find a suitable function $H_t$.  Consider the
  $1$--parameter family of Reeb fields $\smash{R_t}$ as a single vector field
  on the manifold $\smash{[0,1]\times \bigl(\R^3\times T^*S\bigr)}$.  Since
  $\smash{R_t}$ is transverse to the submanifold $\smash{N := [0,1]\times
  \bigl(\R^2\times\{0\}\times T^*S\bigr)}$ along $\smash{[0,1]\times E}$, it is
  possible to define a solution $H_t$ to $dH_t(R_t) = -\dot
  \alpha_t(R_t)$, such that $\left.H_t\right|_{N}\equiv 0$.  In fact,
  because $\left.\dot\alpha\right|_E = 0$, it follows that
  $\left.dH_t\right|_{E} = 0$, and so the vector field $X_t = H_t\,R_t
  + Y_t$ vanishes on $E$.  Hence $X_t$ can be integrated on a small
  neighborhood of $E$, and $E$ is not moved under the flow, which
  finishes the proof of the proposition.
\end{proof}

We can easily choose a compatible almost complex structure $J$ on the
symplectization
\begin{equation*}
  \Bigl(W=\R\times (\Disk{3}\times T^*S),\, \omega =
  d\bigl(e^t\,(dz+\frac{1}{2}\,(x\,dy - y\,dx) + \lcan)\bigr)\Bigr),
\end{equation*}
by observing that the Reeb field is given by $\Reeb =
e^{-t}\partial_z$, and that the kernel of $\alpha$ is spanned by
$\partial_x +\frac{y}{2}\,\partial_z$, $\partial_y
-\frac{x}{2}\,\partial_z$, and the vectors $X-\lcan(X)\,\partial_z$
for all $X\in T(T^*S)$.  Choose a metric $g$ on $S$, and let $J_0$ be
the almost complex structure on $T^*S$ constructed in
\fullref{sec: cotangent bundle} that is compatible with~$d\lcan$.

With this, we can define a $J$ on $W$ by $J\partial_t = \Reeb$,
$J\Reeb = -\partial_t$, $J (\partial_x +\frac{y}{2}\,\partial_z) =
\partial_y -\frac{x}{2}\,\partial_z$, $J (\partial_y
-\frac{x}{2}\,\partial_z) = -\partial_x -\frac{y}{2}\,\partial_z$, and
$J (X -\lcan(X)\,\partial_z) = J_0 X - \lcan(J_0X)\,\partial_z$.  The
last equation can also be written as $JX = J_0 X -
e^t\,\lcan(X)\,\partial_t - \lcan(J_0X)\,\partial_z$.

As a matrix, the complex structure $J$ takes the form:
\begin{align*}
  J(t;x,y,z;\POS,\MOM) &= 
    \begin{pmatrix}
      0  & \frac{y}{2}\,e^t & -\frac{x}{2}\,e^t & -e^t & - e^t\,\lcan \\
      0 & 0 & -1 & 0 & 0 \\
      0 & 1 & 0 & 0 & 0 \\
      e^{-t} & -\frac{x}{2} & -\frac{y}{2} &  0 & -\lcan\circ J_0 \\
      0 & 0 & 0 & 0 & J_0
    \end{pmatrix}
\end{align*}
Note that the last row and column represent linear maps from or to
$T(T^*S)$.  A lengthy computation (which becomes very easy on the
singular set $S$) shows that this structure is compatible
with~$\omega$.

\begin{propo}
  The almost complex manifold $(W,J)$ can be mapped with a
  biholomorphism to
  \begin{align*}
    (\C^2\times T^*S, i\oplus J_0).
  \end{align*}
\end{propo}
\begin{proof}
  The desired biholomorphism is
  \begin{align*}
    \Phi(t,x,y,z;\POS,\MOM) &= (\tilde t,\tilde x,\tilde y,\tilde z;
    \tilde\POS,\tilde\MOM) = \bigg(-e^{-t}-\frac{x^2+y^2}{4}- F,x,y,z;
      \POS,\MOM\bigg),
  \end{align*}
  with the function
  \begin{equation*}
    F\co T^*M \to \R,\quad (\POS,\MOM) \mapsto
    \frac{\norm{\MOM}^2}{2}.
  \end{equation*}
  It brings $J$ into standard form with respect to the coordinate
  pairs $(\tilde x,\tilde y)$, $(\tilde t,\tilde z)$.  More
  explicitly, by pulling back $J$ under the inverse of~$\Phi$
  \begin{align*}
    \Phi^{-1}(\tilde t,\tilde x,\tilde y,\tilde
    z;\tilde\POS,\tilde\MOM) &= (t,x,y,z; \POS,\MOM) =
    \bigg(-\ln(-\tilde t-\frac{\tilde x^2+\tilde y^2}{4}- F),\tilde
      x,\tilde y,\tilde z;\tilde\POS,\tilde\MOM \bigg),
  \end{align*}
  ie\ by computing $D\Phi\cdot J\cdot D\Phi^{-1}$, we obtain the
  matrix
  \begin{align*}
    D\Phi\cdot J\cdot D\Phi^{-1} &=
    \begin{pmatrix}
       0 & 0 & 0 & -1 & -\lcan - dF\circ J_0 \\
       0 & 0 & -1 & 0 & 0 \\
       0 & 1 & 0 & 0 & 0\\  
       1 & 0 & 0 & 0 & dF -\lcan\circ J_0 \\
       0 & 0 & 0 & 0 & J_0
    \end{pmatrix},
  \end{align*}
  and since, according to \fullref{sec: cotangent bundle},
  $dF\circ J_0 = -\lcan$, this gives the desired normal form.
\end{proof}

\begin{figure}[ht!]
\centering
    \begin{picture}(0,0)%
      \includegraphics{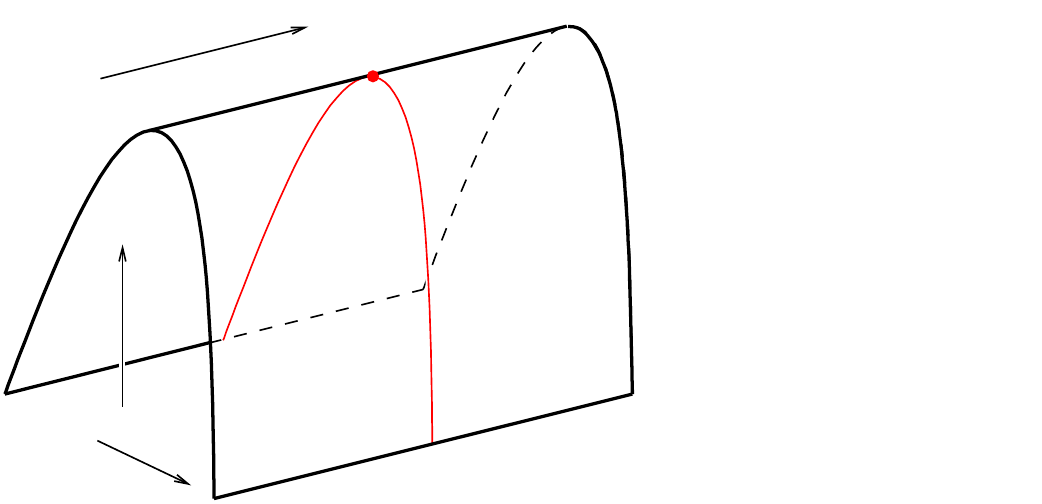} 
    \end{picture}%
    \setlength{\unitlength}{4144sp}%
    \begingroup\makeatletter\ifx\SetFigFont\undefined%
    \gdef\SetFigFont#1#2#3#4#5{%
      \reset@font\fontsize{#1}{#2pt}%
      \fontfamily{#3}\fontseries{#4}\fontshape{#5}%
      \selectfont}%
    \fi\endgroup%
    \begin{picture}(2913,2289)(2208,-4379)
      \put(3800,-2378){\makebox(0,0)[lb]{\smash{{\SetFigFont{6}{7.2}{\familydefault}{\mddefault}{\updefault}{\color[rgb]{1,0,0}$S$}%
            }}}}
      \put(4225,-3984){\makebox(0,0)[lb]{\smash{{\SetFigFont{6}{7.2}{\familydefault}{\mddefault}{\updefault}{\color[rgb]{1,0,0}$F$}%
            }}}}
      \put(3232,-2174){\makebox(0,0)[lb]{\smash{{\SetFigFont{6}{7.2}{\familydefault}{\mddefault}{\updefault}{\adjustlabel<-7pt,0pt>{\color[rgb]{0,0,0}$\ImaginaryPart w_1$}%
            }}}}}
      \put(4726,-4111){\makebox(0,0)[lb]{\smash{{\SetFigFont{8}{9.6}{\familydefault}{\mddefault}{\updefault}{\adjustlabel<-10pt,50pt>{\color[rgb]{0,0,0}$M$}%
            }}}}}
      \put(2791,-3404){\makebox(0,0)[lb]{\smash{{\SetFigFont{6}{7.2}{\familydefault}{\mddefault}{\updefault}{\adjustlabel<-7pt,15pt>{\color[rgb]{0,0,0}$\RealPart w_1$}%
            }}}}}
      \put(2566,-4336){\makebox(0,0)[lb]{\smash{{\SetFigFont{6}{7.2}{\familydefault}{\mddefault}{\updefault}{\adjustlabel<-10pt,0pt>{\color[rgb]{0,0,0}$(w_2,\MOM)$}%
            }}}}}
    \end{picture}%
    \caption{The contact manifold $M=\{\RealPart w_1 = -\frac{1}{2}\,\abs{w_2}^2 -\frac{1}{4}\,\norm{\MOM}^2\}$ appears in the standard neighborhood
      of $S$ like a parabola curved downwards in half of the
      directions.}\label{bild: kontakt mfkt nach unten gekruemmt}
\end{figure}
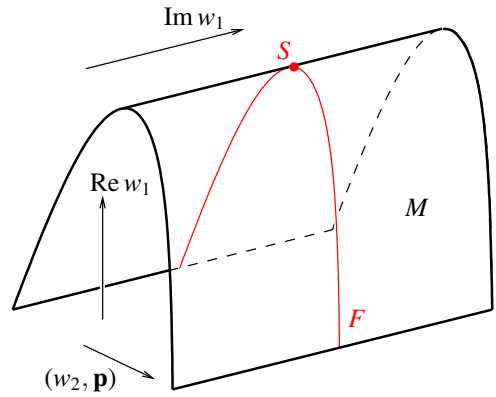

As just proved, the neighborhood of the singular set $S$ in $W$ can be
regarded as $\C^2\times T^*S = \{(w_1,w_2;\POS,\MOM)\}$ with the
almost complex structure $i\oplus J_0$.  The contact manifold $M$ is
given in this model by the set
\begin{equation*}
  \Bigl\{(w_1,w_2;\POS,\MOM)\in \C^2\times T^*S\Bigm|
  \RealPart w_1=-\frac{1}{4}\,\abs{w_2}^2-\frac{1}{2}\,\norm{\MOM}^2\Bigr\}.
\end{equation*}
The contact manifold is thus a hypersurface, represented by a graph,
which is curved downward by the distance from the singular set
(\fullref{bild: kontakt mfkt nach unten gekruemmt}).  By using the
form of a maximally foliated submanifold $F$ found in
\fullref{standard kontaktform bei singularitaet}, one can
write $F$ in a neighborhood of its singular set as the graph
\begin{align*}
  \C\times S & \to \C^2 \times T^*S \\
  (w,\POS) & \mapsto \Bigl(-\frac{1}{4}\,\abs{w}^2,w;\POS,0\Bigr).
\end{align*}

\subsection{Explicit solutions for the Cauchy--Riemann operator}

With the local form obtained in the section above for elliptic
singularities, it is very easy to write down explicitly a family of
holomorphic disks with boundary on the maximally foliated submanifold
$F$.  We will then show that no other simple $J$--holomorphic curves,
which can be capped off with a disk that lies in $F$, can enter this
neighborhood.

Let in this section denote $(W,\omega)$ again a $2n$--dimensional
symplectic filling of $(M,\alpha)$, and let $U$ be the neighborhood of
$S$
\begin{align*}
  \begin{split}
    U = \Bigl\{(z_1,z_2;\POS,\MOM)\in \C^2\times T^*S \Bigm|& -C <
    \RealPart{z_1} \le 0, \quad -C < \ImaginaryPart{z_1} < C, \\
    &\qquad \RealPart{z_1} + \frac{1}{4}\,\abs{z_2}^2 +
    \frac{1}{2}\,\norm{\MOM}^2 \le 0 \Bigr\} \subset W
  \end{split}
\end{align*}
for $C>0$ small enough ($U$ is the half-space below the parabola in
\fullref{bild: trennende mfkt um singularitaet}).

\begin{propo}\label{explizite loesungen} Let us consider the set of
  $J$--holomorphic curves $f\co\Sigma\to W$ whose boundary sits on the
  maximally foliated submanifold $F$.  Among those, the disks
  \begin{align*}
    u_{t_0,\POS_0}\co & \Disk{2} \to \C^2\times T^*S\\
    & z\mapsto (-t_0,2\sqrt{t_0}\,z;\POS_0,0),
  \end{align*}
  for fixed $\POS_0\in S$, and $t_0\in\R_{>0}$ are (up to
  reparametrization) the only simple curves that are completely
  contained in the neighborhood $U$ of the singular set $S$ just
  defined.
\end{propo}
\begin{proof}
  That the maps $u_{t_0,\POS_0}$ are really $J$--holomorphic disks is
  obvious, and that they sit on $F$ can also be checked very easily.

  Assume now, there was a holomorphic curve $u\co\Sigma\to W$
  different from any of the solutions $u_{t_0,\POS_0}$.  If $\partial
  \Sigma = \emptyset$, then $u$ has to be constant, because $\omega$
  is exact in the considered neighborhood.  Note that the projections
  $\pi_1$ and $\pi_2$ of $\C^2\times T^*S$ onto $\C^2$ or onto $T^*S$
  can be concatenated with the map $u$ to provide easier holomorphic
  curves.  In particular, it follows that $u_2 := \pi_2\circ
  u\co\Sigma\to T^*S$ is a $J_0$--holomorphic curve, whose boundary sits
  on the zero-section of the cotangent bundle.  The energy of $u_2$ is
  given by
  \begin{align*}
    E(u_2) &= \int_{u_2} d\lcan = \int_{\Sigma} u_2^* d\lcan =
    \int_{\partial\Sigma} u_2^* \lcan = 0.
  \end{align*}
  Hence it follows that $u_2 = \pi_2\circ u\equiv \POS_0$ is a constant
  curve.

  The other component $\pi_1\circ u$ of the holomorphic curve can be
  written as two ordinary holomorphic functions
  \begin{align*}
    (f_1,f_2) &:= \pi_1\circ u.
  \end{align*}
  The first function $f_1\co\Sigma\to\C$ has vanishing imaginary part
  on all boundary components, and hence it has vanishing imaginary
  part everywhere.  As a consequence, it follows that the real part of
  $f_1$ is constant, as can be seen by using the Cauchy--Riemann
  equations, and deleting the vanishing imaginary part.

  So far, we have shown that any solution is of the form $\smash{z \mapsto
  (-t_0,f_2(z);\POS_0,0)}$, where $f_2\co\Sigma\to\C$ is a holomorphic
  function, such that $\smash{\abs{f_2\big|_{\partial\Sigma}
    }^2 \equiv 4t_0}$.  Assume for simplicity that $t_0 = 1/4$,
  then with the maximum principle it follows that $f_2(\Sigma)
  \subset\Disk{2}$.

  For any $p\in\mathrm{int}\Sigma$, there exist holomorphic charts
  around $p$ and $f_2(p)$ such that $f_2$ takes the form $w\mapsto
  w^k$, but if $k\ne 1$, then $f_2$ is locally a branched covering,
  and $u$ cannot be simple.  It follows that $f_2$ does not have any
  critical points, and so it is a local injective diffeomorphism in
  the interior of $\Sigma$.  In particular, $\Sigma$ has only one
  boundary component, and has to be a disk $\Disk{2}$.
  \begin{align*}
  \tag*{\hbox{The map}}
    h_w\co & \Disk{2} \to \Disk{2} \\
    & z \mapsto \frac{z-w}{1-\bar w z}
  \end{align*}
  is a biholomorphism on $\Disk{2}$, that maps $w\in\Disk{2}$ to $0$.
  If we set $w = f_2(0)$, then the concatenation $H_1:=h_w\circ f_2$
  is a holomorphic map from the unit disk to itself, such that $H_1(0)
  = 0$.  The winding number of $\left.H_1\right|_{\S^1}$ is still $1$.

  Define now a function $H_2(z) := H_1(z)/z$, this map is holomorphic
  on $\Disk{2}-\{0\}$, and can be continuously extended to $0$ by
  setting $\smash{H_2(0) := H_1^\prime(0)}$.  This extension is also
  holomorphic. Furthermore, $\smash{H_2(\S^1)\subset \S^1}$, and by the
  maximum principle $H_2(\Disk{2}) \subset \Disk{2}$.  The winding
  number of $\left.H_2\right|_{\S^1}$ is zero, and hence it follows
  that $\partial_\phi H_2(e^{i\phi})$ vanishes for some angle
  $\phi_0$.  With the Cauchy--Riemann equation it follows that
  $\partial_r H_2(e^{i\phi}) = 0$, because
  \begin{align*}
    \partial_\phi H_2 &= (x\,\partial_y - y\,\partial_x)\,(\RealPart
    H_2 + i\ImaginaryPart H_2) = i\,\partial_r H_2.
  \end{align*}
  In particular it follows that the function $\abs{H_2}^2$ has
  vanishing derivative at $e^{i\phi_0}$, but this contradicts the
  boundary point lemma \cite[Lemma~3.4]{GilbargTrudinger}, and
  hence $H_2\equiv\mathrm{const}$ and $H_1(z) = e^{i\theta_0}z$.
\end{proof}

The next remark allows us to apply \fullref{kurven bleiben weg
  von der singularitaet} in our situation.

\begin{remark}
  Any $J$--holomorphic disk $u$ that lies in the moduli space of
  curves starting at a Bishop family, can be capped off by a
  (topological) disk that is completely contained in the maximally
  foliated submanifold~$F$.
\end{remark}

\begin{propo}\label{kurven bleiben weg von der singularitaet} The only
  nontrivial holomorphic curves intersecting the neighborhood~$U$
  that can be capped off by attaching disks lying in the maximally
  foliated submanifold $F$ are the curves $u_{t_0,\POS_0}\!$ (and their
  multiple covers) defined in \hbox{\fullref{explizite loesungen}}.
\end{propo}
\begin{proof}
  Suppose first that the restriction $\smash{\left.u\right|_{u^{-1}(U)}}$ of
  the curve $u$ has constant $(\ImaginaryPart z_1)$--coordinate on one
  component of $\smash{u^{-1}(U)}$.  Then the whole curve $u$ is contained in
  $U$, because with the Cauchy--Riemann equation, it follows already
  that $z_1$ is constant on this component of $u^{-1}(U)$.  If $u$
  approaches the boundary of the neighborhood $U$, then either the
  $\MOM$-- or the $z_2$--coordinate have to grow, but as soon as
  $\RealPart z_1 + \frac{1}{4}\,\abs{z_2}^2 +
  \frac{1}{2}\,\norm{\MOM}^2$ vanishes, $u$ touches the contact
  manifold $M$, so that $u$ is trapped by the hypersurface $M$ inside
  $U$.

  By \fullref{explizite loesungen}, this means that the only
  nontrivial holomorphic curves~$u$ intersecting the neighborhood
  $U$, having constant $(\ImaginaryPart z_1)$--coordinate on a
  component of $u^{-1}(U)$, are the disks $u_{t_0,\POS_0}$ given above
  (and their multiple covers).
  
  Assume now that a curve $u$ enters $U$, but does not have constant
  $(\ImaginaryPart z_1)$--coordinate.  To disprove the existence of
  $u$, we will use an intersection argument similar to the classical
  one in dimension~$4$.
  
  Consider for every $c=c_x+ic_y\in\C$ the submanifold (as drawn in
  \fullref{bild: trennende mfkt um singularitaet})
  \begin{align*}
    A_c &:= \Bigl( \bigl\{x+ic_y\in\C\bigm| x\ge c_x\}\times\C\times
    T^*S \Bigr) \cap U.
  \end{align*}
  Such an $A_c$ is the codimension~$1$ submanifold in~$U$ obtained by
  taking a slice with constant imaginary $z_1$--coordinate and
  chopping off everything having smaller real $z_1$--coordinate
  than~$c_x$.  The boundary of $A_c$ is composed by two smooth
  manifolds: one of them is
  \begin{align*}
    B_c &:= \Bigl( \bigl\{c\}\times\C\times T^*S \Bigr) \cap U
  \end{align*}
  which is a compact $J$--holomorphic codimension~$2$ submanifold.
  Note that $U$ is foliated by the $B_c$ for different values of
  $c\in\C$.

  The second part of the boundary is given by the set $A_c\cap M$.
  Since the boundary is convex, holomorphic curves can only touch
  $A_c\cap M$ at their own boundary, but since the boundary of the
  holomorphic curves $u$ we are considering, lies in $F$, $u$ will
  never intersect $A_c\cap M$ if $c_y\ne 0$.

  Thus, we have obtained a (nonsmooth) closed manifold $\widehat B_c
  =\partial A_c$, that represents the trivial homology class in
  $H_{2n-2}(W)$, and which allows us to compute the intersection
  number with holomorphic curves.

  \begin{figure}[ht!]
    \begin{center}
      \begin{picture}(0,0)%
        \includegraphics{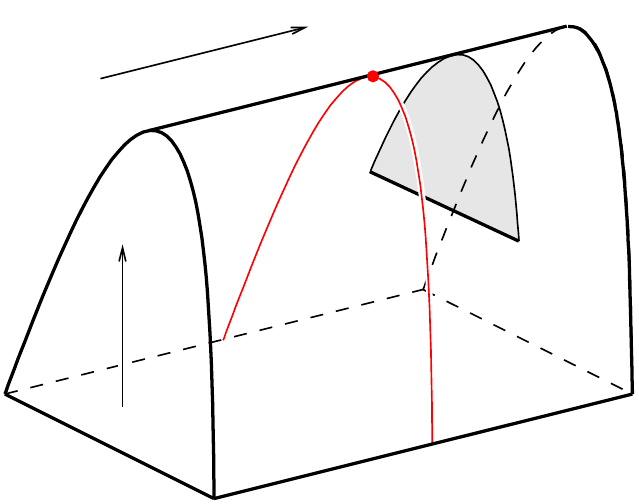}%
      \end{picture}%
      \setlength{\unitlength}{4144sp}%
      \begingroup\makeatletter\ifx\SetFigFont\undefined%
      \gdef\SetFigFont#1#2#3#4#5{%
        \reset@font\fontsize{#1}{#2pt}%
        \fontfamily{#3}\fontseries{#4}\fontshape{#5}%
        \selectfont}%
      \fi\endgroup%
      \begin{picture}(2913,2289)(2208,-4379)
        \put(3800,-2378){\makebox(0,0)[lb]{\smash{{\SetFigFont{6}{7.2}{\familydefault}{\mddefault}{\updefault}{\color[rgb]{1,0,0}$S$}%
              }}}}
        \put(4225,-3984){\makebox(0,0)[lb]{\smash{{\SetFigFont{6}{7.2}{\familydefault}{\mddefault}{\updefault}{\color[rgb]{1,0,0}$F$}%
              }}}}
        \put(3232,-2174){\makebox(0,0)[lb]{\smash{{\SetFigFont{6}{7.2}{\familydefault}{\mddefault}{\updefault}{\adjustlabel<-4pt,0pt>{\color[rgb]{0,0,0}$\ImaginaryPart z_1$}}%
              }}}}
        \put(3331,-3976){\makebox(0,0)[lb]{\smash{{\SetFigFont{8}{9.6}{\familydefault}{\mddefault}{\updefault}{\color[rgb]{0,0,0}$M$}%
              }}}}
        \put(2791,-3404){\makebox(0,0)[lb]{\smash{{\SetFigFont{6}{7.2}{\familydefault}{\mddefault}{\updefault}{\color[rgb]{0,0,0}$\RealPart z_1$}%
              }}}}
        \put(4286,-3186){\makebox(0,0)[lb]{\smash{{\SetFigFont{6}{7.2}{\familydefault}{\mddefault}{\updefault}{\adjustlabel<0pt,-4pt>{\color[rgb]{0,0,0}$B_c$}}%
              }}}}
        \put(4214,-2755){\makebox(0,0)[lb]{\smash{{\SetFigFont{6}{7.2}{\familydefault}{\mddefault}{\updefault}{\adjustlabel<-5pt,0pt>{\color[rgb]{0,0,0}\colorbox{white}{$A_c$}}%
              }}}}}
      \end{picture}%
      \caption{The neighborhood $U$ of the singular set is foliated by
        $J$--holomorphic codimension~$2$ submanifolds $B_c$, which
        represent one part of the boundary of the slices $A_c$.  This
        gives rise to an intersection argument.}\label{bild: trennende
        mfkt um singularitaet}
    \end{center}
  \end{figure}
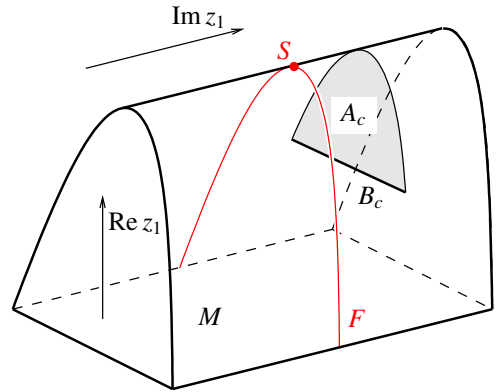

  Let $u$ be now any holomorphic curve that passes through the model
  neighborhood $U$, and whose boundary lies on the maximally foliated
  submanifold~$F$.  Assume $u$ can be capped off by attaching a disk
  in~$F$, and denote $u$ together with its attached disk by~$\widehat
  u$.

  If $u$ has nonconstant imaginary $z_2$--coordinate in the model
  neighborhood $U$, then there is (by Sard's Theorem) a submanifold
  $B_c$ (with $\ImaginaryPart c \ne 0$) that intersects $u$
  transversally at a discrete set of points.  The two homology classes
  $[\widehat u]$ and $[\widehat B_c]$ have positive intersection
  number, because the only intersections between both classes lie in
  the subset, where both are represented by $J$--holomorphic
  submanifolds.  But a positive intersection number is not possible,
  since $[\widehat B_c]$ represents the trivial homology class in
  $H_{2n-2}(W)$.
\end{proof}

\subsection{Expected dimension for the bishop family}

The expected dimension for the solution space of the Cauchy--Riemann
operator at a holomorphic disk $u$, whose boundary lies on a totally
real submanifold $F$, is given by the formula (see
McDuff and Salamon \cite[Theorem~C.1.10]{McDuffSalamonJHolo})
\begin{align*}
  \mathrm{index}\, \bar\partial_J &= \frac{\dim W}{2} + \mu(u^*TW,
  u^*TF),
\end{align*}
where we have used that the Euler characteristic of a disk is
$\chi(\Disk{2}) = 1$.  Here $\mu(E_\C,E_\R)$ denotes the Maslov index
of a complex vector bundle $E_\C$ over a disk $\Disk{2}$ with respect
to a totally real subbundle $E_\R\le \left.E_\C\right|_{\S^1}$ over
the boundary of~$\Disk{2}$.

To obtain the expected dimension of the moduli space, we have to
subtract $2$, which corresponds to the dimension of the automorphism
group of the holomorphic unit disk with one marked point on the
boundary.

\begin{propo}\label{maslov index der bishop kurven}
  The Maslov index $\mu(u^*TW, u^*TF)$ is equal to $2$ for any of the
  holomorphic disks
  \begin{align*}
    u\co & \Disk{2} \hookrightarrow \C^2\times T^*S\\
    & z\mapsto (-t_0,2\sqrt{t_0}\,z;\POS_0,0)
  \end{align*}
  given in \fullref{explizite loesungen} above.
\end{propo}
\begin{proof}
  We can trivialize $u^*TW$ by choosing the obvious complex basis
  \begin{align*}
    u^*TW &= \langle \partial_{x_1},\partial_{x_2},\partial_{q^1},
    \dotsc, \partial_{q^n}\rangle_\C.
  \end{align*}
  The totally real subbundle at the boundary of $u(e^{i\phi})$ is
  spanned by the vectors
  \begin{align*}
    \left(\left.u\right|_{\partial\Disk{2}}\right)^*TF &= \langle
    \partial_{x_2} - \sqrt{t_0}\cos\phi\,\partial_{x_1},\partial_{y_2}
    - \sqrt{t_0}\sin\phi\,\partial_{x_1}, \partial_{q^1}, \dotsc,
    \partial_{q^n} \rangle_\R.
  \end{align*}
  This subbundle can be represented at a point
  $(-t_0,2\sqrt{t_0}\,e^{i\phi};\POS_0)\in F$ by the matrix
  \begin{align*}
    \Lambda(e^{i\phi}) &=
      \begin{pmatrix}
        -\sqrt{t_0}\cos\phi & -\sqrt{t_0}\sin\phi & 0 \\
        1 & i & 0 \\
        0 & 0 & \1
      \end{pmatrix}
  \end{align*}
  with respect to the complex basis of $u^*TW$.  It follows that the
  Maslov index $\mu(u^*TW, u^*TF)$ is given by
  $$\eqalignbot{
    \mu(u^*TW, u^*TF) &= \deg \frac{\det \Lambda^2}{\det
      \Lambda^*\Lambda} =
    \deg \frac{-t_0\,(\cos^2\phi - \sin^2\phi + 2i\sin\phi\cos\phi)}{t_0}\cr
    &= \deg (-e^{2i\phi}) = 2. \cr
  } \proved $$
\end{proof}

The expected dimension of the moduli space $\mathcal{M} =
\mathcal{M}(\C^2\times T^*S,z_0)$ is thus
\begin{align*}
  \dim\mathcal{M} &= \frac{\dim W}{2} + 2 - 2 = \frac{\dim W}{2},
\end{align*}
which means that the family of solutions, we have found above are
locally all solutions, if the Cauchy--Riemann operator for the given $J$
is regular.  That this is indeed the case will be shown below.

\subsection{Surjectivity of the linearized Cauchy--Riemann operator}
\label{sec: regularitaet von CR in naehe von singularitaet}

In this section, we will show that the linearized Cauchy--Riemann
operator $\bar\partial_J$ at any of the curves $u_{t_0,\POS_0}$ from
\fullref{explizite loesungen} is surjective.  Since we have
shown that these disks are the only holomorphic curves contained in
our model neighborhood, we can apply
\cite[Remark~3.2.3]{McDuffSalamonJHolo}.  Thus it will be enough to
perturb $J$ outside the model neighborhood to obtain a regular
Cauchy--Riemann operator for the whole moduli space.

\begin{propo}
  The linearized Cauchy--Riemann operator $D_u\bar\partial_J$ is
  surjective at any of the disks $u = u_{t_0,\POS_0}$ specified in
  \fullref{explizite loesungen}.
\end{propo}
\begin{proof}
  Instead of checking that $D_u\bar\partial_J$ is surjective, we will
  compute the dimension of its kernel, and see that it coincides with
  the index of the operator.  It then follows that the cokernel must
  be trivial.

  Write again $W$ for the neighborhood $\C^2\times T^*S$, and let
  $v+iw$ be the standard coordinates on the disk.  The linearized
  Cauchy--Riemann operator at a disk $u\co\Disk{2}\to W$ is given by
  \begin{align*}
    (D_u\bar \partial_J) \dot u &= \frac{\partial \dot u}{\partial v}
    + J(u)\,\frac{\partial\dot u}{\partial w} +
    \left(\left.\frac{d}{ds}\right|_{s=0} J(u+s\dot
      u)\right)\,\frac{\partial u}{\partial w},
  \end{align*}
  where $\dot u$ denotes a section of $u^*TW$, which restricts to a
  section of $u^* TF$ along the boundary of~$\Disk{2}$.  To compute
  the kernel, we have to find all solutions $\dot u$ of the equation
  $(D_u\bar \partial_J)\dot u = 0$.

  Note that the linearized Cauchy--Riemann equation simplifies to
  \begin{align*}
    \frac{\partial \dot u}{\partial v} + J(u)\,\frac{\partial\dot
      u}{\partial w} &= 0,
  \end{align*}
  because the derivative of the almost complex structure can be
  dropped: for the $\C^2$ part, this is obvious, because the complex
  structure is constant; for the $T^*S$ part, this follows, because
  we multiply with $\partial_w u$, which vanishes in the cotangent
  bundle for any of the curves $u = u_{t_0,\POS_0}$ from
  \fullref{explizite loesungen}.

  Take a canonical chart $\{(z_1,z_2;\POS,\MOM)\}$ in $W$ containing
  the disk $u$, and use this to express $\dot u$ as
  \begin{align*}
    \dot u\co & \Disk{2} \to \C^2\times\R^{2n-4} \\
    & v+iw \mapsto (\dot z_1,\dot z_2; \dot\POS,\dot\MOM).
  \end{align*}
  The boundary conditions give 
  $\left.\ImaginaryPart(\dot z_1)\right|_{\S^1} \equiv 0$,
  $\left.\RealPart (\dot z_1)\right|_{\S^1} =
  \frac{1}{4}\,\big(z_2\bar{\dot z}_2 + \bar z_2 \dot z_2\big)$
  and $\left.\dot\MOM\right|_{\S^1} \equiv 0$.

  The linearized Cauchy--Riemann equation decomposes into two
  independent equations, one on~$\C^2$, and one on~$T^*S$.  We will
  first analyze the part on~$T^*S$.  The holomorphic disk at which we
  are linearizing is just the constant map $(\POS_0,0)\in T^*S$.  We
  can assume the chart $(q^1,\dotsc,q^{n-2},p_1,\dotsc,p_{n-2})$ to be
  induced by a geodesic normal chart $(q^1,\dotsc,q^{n-2})$ around
  $\POS_0\in S$.  The explicit form of $J_0$ in such a chart can be
  found in the proof of \fullref{laenge auf faser J und lcan} in
  \fullref{sec: cotangent bundle}.  The linearized Cauchy--Riemann
  equation becomes
  \begin{equation*}
    \frac{\partial \dot\POS}{\partial v} - \frac{\partial\dot\MOM}{\partial w} = 0
    \qquad \text{ and } \qquad
    \frac{\partial \dot\MOM}{\partial v} + \frac{\partial\dot\POS}{\partial w} = 0.
  \end{equation*}
  This gives rise to the equation $\Delta\dot\MOM = 0$, which in turn
  together with the boundary condition implies that $\dot\MOM \equiv
  0$.  By plugging this into the Cauchy--Riemann equation, it finally
  follows that $\dot\POS \equiv \mathrm{const}$.  That means we find
  $n-2$ degrees of freedom from the $T^*S$ part.

  The equation for  $z_1$ is the standard Cauchy--Riemann equation, ie\ 
  \begin{align*}
    \frac{\partial \dot z_1}{\partial v} + i \,\frac{\partial\dot
      z_1}{\partial w} &= 0,
  \end{align*}
  or split into real part and imaginary part (with $\dot z_1 = \dot
  x_1 + i\dot y_1$)
  \begin{equation*}
    \frac{\partial \dot x_1}{\partial v} -\frac{\partial\dot
      y_1}{\partial w} = 0\quad \text{ and }\quad
    \frac{\partial \dot y_1}{\partial v} + \frac{\partial\dot
      x_1}{\partial w} = 0.
  \end{equation*}
  Combining these equations, we obtain $\Delta \dot y_1 \equiv 0$, and
  together with the boundary condition, it follows that $\dot y_1
  \equiv 0$.  Using this result again in the linearized Cauchy--Riemann
  equation, we obtain $\dot x_1 \equiv \mathrm{const}$.

  The function $\dot z_2$ is a holomorphic function on the unit disk,
  ie\ we can write $\dot z_2$ down as a power series
  \begin{align*}
    \dot z_2(v+iw) :&= \sum_{k=0}^\infty a_k\, (v+iw)^k,
  \end{align*}
  or when restricted to the boundary $\S^1 = \{e^{i\phi}\}$ and
  assuming that $a_k = b_k + ic_k$,
  \begin{align*}
    \dot z_2(e^{i\phi}) :&= \sum_{k=0}^\infty (b_k+ic_k)\, e^{ik\phi}.
  \end{align*}
  The function $\dot z_1 \equiv\mathrm{const}$ was coupled to $\dot
  z_2$ by the boundary condition
  \begin{align*}
    \left.\dot z_1\right|_{\S^1} &=
    \frac{1}{4}\,\left.\Bigl(z_2\bar{\dot z}_2 + \bar z_2 \dot
      z_2\Bigr)\right|_{\S^1},
  \end{align*}
  and by using that $z_2 = \sqrt{t_0}\,(v+iw)$ (after a
  reparametrization), it follows that
  \begin{align*}
    \mathrm{const} &\equiv e^{-i\phi}\,\sum_{k=0}^\infty (b_k+ic_k)\,
    e^{ik\phi} + e^{i\phi}\,\sum_{k=0}^\infty (b_k-ic_k)\,
    e^{-ik\phi}\\
    &= \sum_{k=0}^\infty (b_k+ic_k)\, e^{i\,(k-1)\,\phi} +
    \sum_{k=0}^\infty (b_k-ic_k)\, e^{-i\,(k-1)\,\phi} \\
    &= \sum_{k=-1}^\infty (b_{k+1}+ic_{k+1})\, e^{ik\phi} +
    \sum_{k=-1}^\infty (b_{k+1}-ic_{k+1})\, e^{-ik\phi}\\
    &= 2 \sum_{k=-1}^\infty b_{k+1}\cos k\phi - 2 \sum_{k=-1}^\infty
    c_{k+1}\sin k\phi \\
    \begin{split}
      &= 2 b_1 - 2c_1 + 2\, (b_0+b_2)\,\cos\phi + 2\,(c_0-c_2)\,\sin\phi \\
      &\hphantom{= 2 b_1\ } + 2\, \sum_{k=2}^\infty \Bigl( b_{k+1}\cos k\phi -
      c_{k+1}\sin k\phi \Bigr).
    \end{split}
  \end{align*}
  And hence the coefficients $a_0$ and $a_1$ can be chosen
  arbitrarily, $a_2 = -\bar a_0$, and $a_3=a_4=\dotsc=0$ must all
  vanish. The function $\dot z_2$ is thus given by $\dot z_2 = a_0 +
  a_1 z - \bar a_0 z^2$, and so the $\C^2$ part contributes a
  $4$--dimensional kernel to the linearized Cauchy--Riemann operator.

  Since both parts are independent, the dimension of the kernel of the
  linearized Cauchy--Riemann operator is equal to $4+\dim S =
  \unfrac{\dim W}{2} + 2$, which is equal to the Fredholm index, as we
  wanted to show.
\end{proof}

\section{Bubbling off analysis}\label{sec: no bubbling}

The moduli spaces of holomorphic disks, whose boundary lies on a
compact totally real submanifold, and which have a uniform energy
bound, are compact, provided one includes two different types of
bubbling: the curves can either form a bubble at their boundary or a
bubble in their interior.  Though a maximally foliated submanifold $F$
with elliptic singularities is not compact, if one removes the
singularities (or not totally real, if one does not remove the
singularities), we still have compactness for the moduli spaces coming
from Bishop families, because we have proved in
\hbox{\fullref{kurven bleiben weg von der singularitaet}} that
elliptic singularities have a small neighborhood, which blocks out
every curve with exception of a Bishop family.  This confines the
movement of all other disks to a compact totally real submanifold.
But first, we will show that the energy for every curve of the moduli
space is bounded.

\begin{propo}\label{energie ist beschraenkt} Let $(M^{2n-1},\alpha)$
  be a closed contact manifold that has a symplectic filling
  $(W,\omega)$.  Assume that $M$ contains a plastikstufe
  $\plastikstufe{S}$.  Let $u\co (\Disk{2},\S^1)\to
  (W,\plastikstufe{S})$ be a holomorphic disk that lies in the same
  moduli space as the Bishop family found in
  \fullref{explizite loesungen}.

  There is a constant $C$, which only depends on
  $\left.\alpha\right|_{T\plastikstufe{S}}$ that bounds the energy
  \begin{align*}
    E(u) &= \int_u \omega
  \end{align*}
  of any such disk.
\end{propo}
\begin{proof}
  The disk $u$ is in $(W,\plastikstufe{S})$ homotopic to a point,
  hence the energy $E(u)$ of $u$ can be obtained by
  \begin{align*}
    E(u) &= \int_u \omega = \int_{\partial u} \alpha.
  \end{align*}
  Let now $P$ denote the submanifold $\plastikstufe{S} - (S \cup
  \partial\plastikstufe{S})$, ie\ the plastikstufe with its boundary
  and the interior singularity removed.  Note that the standard leaves
  of the foliation of $\plastikstufe{S}$ can be labeled bijectively
  with elements $e^{i\phi}\in\S^1$, and in fact there is a smooth
  surjective map
  \begin{align*}
    \theta\co & P
    \to \S^1,
  \end{align*}
  such that $\theta(p) = \theta(q)$ if and only if, $p$ and $q$ lie on
  the same leaf of the foliation.  The differential $d\theta$ can be
  regarded as a $1$--form on $P$.  Note that there is a unique smooth
  function $f\co P\to \R_{>0}$ such that
  \begin{align*}
    \alpha &= f\,d\theta,
  \end{align*}
  because the kernel of both $1$--forms agree on $P$.

  The boundary of every holomorphic disk $u$ in our moduli space is
  transverse to the foliation, because if the tangent direction of $u$
  at a boundary point $z\in\partial\Disk{2}$ lay in $\ker\alpha$,
  then by our definition of $J$ the whole tangent space $T_{u(z)} u$
  would be tangent to $M$.  But this contradicts \fullref{kurven
    transvers zu rand}, and hence it follows that $\theta\circ
  \left.u\right|_{\S^1}$ makes exactly one turn, or expressed in a
  different way,
  \begin{align*}
    \int_{\partial u} d\theta = 2\pi.
  \end{align*}
  For every point $p_0$ on $\partial\plastikstufe{S}$, we find a chart
  $U\subset \R_{\le 0}\times \R^{N-1}$ ($N:=\dim\plastikstufe{S}$)
  mapping $p_0$ to $0$ such that the intersection of the leaves with
  $U$ is given by the planes
  \begin{align*}
    F_{x} &= \Bigl(\{x\}\times \R^{N-1}\Bigr) \cap U,
  \end{align*}
  for $x\le 0$, and $F_0$ corresponds to the boundary
  $\partial\plastikstufe{S}\cap U$.  The local picture alone does not
  allow us to see, which two sheets are contained in the same leaf
  (see \hbox{\fullref{bild: blaetterung am rand der plastikstufe}}), but
  since the leaves in $P$ approach the boundary of the plastikstufe
  asymptotically, it follows that for every leaf, there is a
  monotonous sequence $a_k$ converging to $0$ such that all $F_{a_k}$
  lie in the same global leaf (and such that every other
  hyperplane~$F_b$ lies in another leaf).
  
  \begin{figure}[ht!]
    \begin{center}
      \begin{picture}(0,0)%
        \includegraphics{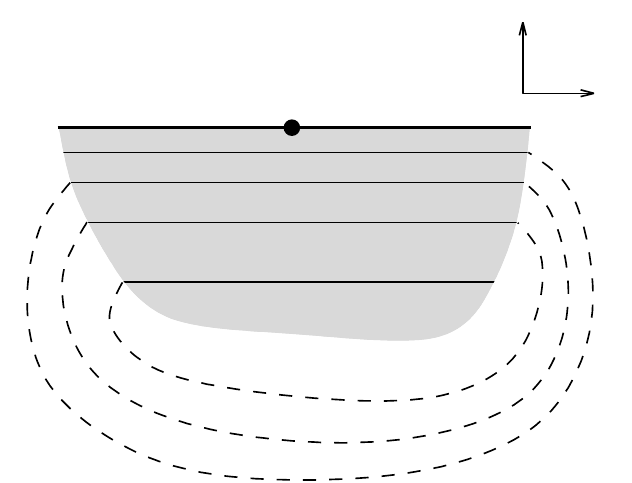}%
      \end{picture}%
      \setlength{\unitlength}{4144sp}%
      \begingroup\makeatletter\ifx\SetFigFont\undefined%
      \gdef\SetFigFont#1#2#3#4#5{%
        \reset@font\fontsize{#1}{#2pt}%
        \fontfamily{#3}\fontseries{#4}\fontshape{#5}%
        \selectfont}%
      \fi\endgroup%
      \begin{picture}(2876,2206)(1913,-3911)
        \put(4189,-1789){\makebox(0,0)[lb]{\smash{{\SetFigFont{8}{9.6}{\familydefault}{\mddefault}{\updefault}{\adjustlabel<0pt,1pt>{\color[rgb]{0,0,0}$\R$}%
              }}}}}
        \put(4668,-2177){\makebox(0,0)[lb]{\smash{{\SetFigFont{8}{9.6}{\familydefault}{\mddefault}{\updefault}{\color[rgb]{0,0,0}$\R^{N-1}$}%
              }}}}
        \put(3206,-2200){\makebox(0,0)[lb]{\smash{{\SetFigFont{8}{9.6}{\familydefault}{\mddefault}{\updefault}{\adjustlabel<0pt,1pt>{\color[rgb]{0,0,0}$p_0$}%
              }}}}}
        \put(3526,-2832){\makebox(0,0)[lb]{\smash{{\SetFigFont{8}{9.6}{\familydefault}{\mddefault}{\updefault}{\adjustlabel<19pt,-1pt>{\color[rgb]{0,0,0}$F_{\smash{x_2}}$}%
              }}}}}
        \put(3509,-2645){\makebox(0,0)[lb]{\smash{{\SetFigFont{8}{9.6}{\familydefault}{\mddefault}{\updefault}{\adjustlabel<38pt,-1pt>{\color[rgb]{0,0,0}$F_{\smash[t]{x_3}}$}%
              }}}}}
        \put(3571,-3108){\makebox(0,0)[lb]{\smash{{\SetFigFont{8}{9.6}{\familydefault}{\mddefault}{\updefault}{\adjustlabel<0pt,-1pt>{\color[rgb]{0,0,0}$F_{x_1}$}%
              }}}}}
      \end{picture}%
      \caption{In a neighborhood of a point $p_0\in\plastikstufe{S}$,
        the leaves can be represented by horizontal hyperplanes $F_x$.
        The global shape of the foliation in $P$ connects certain
        $F_x$.}\label{bild: blaetterung am rand der plastikstufe}
    \end{center}
  \end{figure}
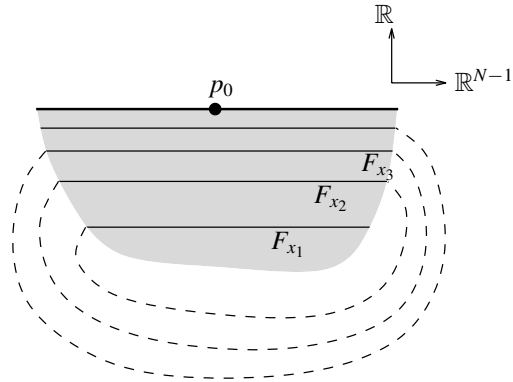
  
  Since $\partial\plastikstufe{S}$ is compact, we can cover the whole
  boundary of the plastikstufe by using only a finite number of such
  charts $\{U_1,\dotsc,U_{n_0}\}$, $n_0\in\N$.  On each chart $U_k$,
  $\alpha$ can be written as $g_k\,dx_1$, where $g_k\co U_k\to\R$ is a
  smooth function.  In particular, $\abs{g_k}$ is bounded by a
  number~$c_k$.

  For any smooth path $\gamma$ in $U_k$ connecting $F_{a}$ with
  $F_{b}$ crossing each leaf transversely in increasing direction,
  we can estimate the integral by
  \begin{equation*}
    \int_{\gamma} \alpha \le c_k\,\abs{b-a}.
  \end{equation*}
  The intersection $\partial u\cap U_k$, where $\partial u$ is the
  boundary of a $J$--holomorphic disk, gives a collection of paths
  $\gamma_j\subset U_k$.  But since $\partial u$ crosses every global
  leaf in the plastikstufe once, we can order the segments of paths in
  such a way that the end point of $\gamma_j$ lies on $F_{b_j}$, and
  the segment $\gamma_{j+1}$ starts on $F_{a_{j+1}}$ with $b_j <
  a_{j+1}$.  Hence the total estimate gives
  \begin{equation*}
    \int_{\partial u\cap U_k} = \sum_j \int_{\gamma_j}\alpha
    \le c_k\,\sum_j\abs{b_j-a_j} < c_kL_k,
  \end{equation*}
  where $L_k$ is the length of $U_k$ in $x_1$--direction.

  Finally cover now the whole plastikstufe with sets
  $\{K,U_1,\dotsc,U_{n_0}\}$, where $K$ is a compact set that never
  touches the boundary of the plastikstufe.

  The energy of a holomorphic disk $u$ in our moduli space with
  boundary $\partial u$ can now be bounded by a number $C$ in the
  following way:
  \begin{align*}
    E(u) &= \int_{\partial u} \alpha \le \int_{\partial u\cap K}
    \alpha + \sum_{k=1}^{n_0} \int_{\partial u\cap U_k} \alpha \le
    \int_{\partial u\cap K} \alpha + \sum_{k=1}^{n_0} c_kL_k.
  \end{align*}
  The energy of the segments of $\partial u$ contained in $K$ can be
  estimated by
  \begin{align*}
    \int_{\partial u\cap K} \alpha &= \int_{\partial u\cap K}
    f\,d\theta \le \sup_{p\in K} f(p)\,\int_{\partial u\cap K} d\theta
    \le 2\pi\, \sup_{p\in K} f(p).
  \end{align*}
  By setting $C:= 2\pi\, \sup_K f + \sum_{k=1}^{n_0} c_kL_k$, we
  obtain a uniform estimate for all disks in our moduli space.
\end{proof}

\begin{propo}
  Let $(M^{2n-1},\alpha)$ be a closed contact manifold that has a
  symplectic filling $(W,\omega)$, and assume that $M$ contains a
  plastikstufe $\plastikstufe{S}$.

  Let $u_k\co\Disk{2}\to W$ be a sequence of holomorphic disks, whose
  boundary lies in the plastikstufe with fixed linking number
  $\linkingzahl(\partial u_k,S) = 1$.  The sequence $u_k$ has a
  subsequence that converges either to a constant map, whose image
  lies in the singular set $S$, or to a simple smooth holomorphic disk
  and a finite number of bubble spheres.
\end{propo}
\begin{proof}
  Assume first, there was a subsequence of disks
  $\bigl(u_{k_l}\bigr)_l$ coming arbitrarily close to the singular
  set.  By \fullref{kurven bleiben weg von der singularitaet}, it
  follows that for $k_l$ sufficiently large, the disks $u_{k_l}$ lie
  in the Bishop family, and then $\bigl(u_{k_l}\bigr)_l$ has a further
  subsequence that converges to a point of $S$.

  If all the disks stay at a finite distance $C>0$ from $S$, then the
  boundary of all disks is contained in a compact subset of a totally
  real submanifold, and hence by Gromov compactness, there is a
  subsequence $\bigl(u_{k_l}\bigr)_l$ that converges to a bubbled
  curve.

  One possible type of bubbles that could occur, are disks growing at
  the boundary of the family $\bigl(u_{k_l}\bigr)_l$.
Let 
  \begin{align*}
    \theta\co & \plastikstufe{S} - (S \cup \partial\plastikstufe{S})
    \to \S^1
  \end{align*}
  be the function already defined in the proof of
  \fullref{energie ist beschraenkt}.  Since the boundary of
  every disk $u_k$ is transverse to the foliation, and the linking
  number is~$1$, we can define smooth bijective maps $\smash{f_k :=
  \theta\circ \left.u_k\right|_{\S^1}\co \partial\Disk{2} \to \S^1}$.
  The sequence $\smash{\bigl(u_{k_l}\bigr)_l}$ converges to a continuous map
  $u_\infty\co \Disk{2}\to W$ in $C^0$--norm, and
  $\bigl(f_{k_l}\bigr)_l$ converges then to $f_\infty := \theta\circ
  \left.u_\infty\right|_{\S^1}$.  The map $f_\infty$ is continuous,
  monotonous, and has degree $1$. It follows that the only way to
  split $u_\infty$ into subbubbles is by assuming that one of the
  bubbles is constant at its boundary, which also implies by the same
  energy argument as above that the bubble is constant.
  
  This means that the only type of bubbling, which is allowed, is
  bubbling in the interior of the disk.  In this situation, all the
  bubbles are spherical.  The base of the bubble tree is a simple disk
  $u_0$, because the restriction to the boundary of~$u_0$ is by the
  argument above injective.
\end{proof}

\begin{defi}
  Let $A$ be a subset of a manifold $M$.  We say that $A$ has
  \textit{at most dimension $n$}, if there is an $n$--dimensional
  manifold $X$ (with finite amount of components) and a smooth map
  $f\co X\to M$ with closed image such that $A\subset f(X)$.
\end{defi}

\begin{remark}
  To apply the general theory of moduli spaces of $J$--holomorphic
  curves, we always assume regularity of the Cauchy--Riemann operator
  for every curve in $\mathcal{M}(W,\plastikstufe{S},z_0)$ and for all
  the bubble trees considered in the proof of the next proposition.
  Below we will briefly sketch the argument, why this is indeed
  possible.

  The almost complex structure $J$ on $W$ is first defined only in a
  neighborhood of the singular set $S$ of $\plastikstufe{S}$ as
  explained in \fullref{sec: regularitaet von CR in naehe von
    singularitaet}.  In this neighborhood, regularity works as proved
  for $\mathcal{M}(W,\plastikstufe{S},z_0)$ and since no holomorphic
  spheres can enter the domain, no problem occurs.

  In a second step $J$ is extended to a small collar neighborhood of
  $M$, such that $J$ is compatible with the convex boundary as defined
  in \fullref{sec: wiederholung von fakten}.  Again holomorphic
  spheres pose no problem, and since all holomorphic disks in our
  moduli space are simple, by perturbing $J$, one can achieve
  regularity for all disks contained in the collar.

  In the last step, $J$ is finally extended over the rest of the
  symplectic manifold $W$, now only requiring that $J$ is compatible
  with~$\omega$.  For disks, one could suspect a difficulty, because
  the boundary part of the disk lies in the collar where $J$ has
  already been defined, but \cite[Remark~3.2.3]{McDuffSalamonJHolo}
  tells us that, regularity for these curves can be achieved even by
  perturbing $J$ only in the interior of~$W$.  For the bubble trees,
  we also obtain by perturbations of $J$ regularity as explained in
  \cite{McDuffSalamonJHolo}.
\end{remark}

\begin{propo}
  To compactify the moduli space
  $\mathcal{M}(W,\plastikstufe{S},z_0)$, one has to add bubbled
  curves.  The image of these bubbled curves under the evaluation map
  $\evaluationmap_{z_0}$ has at most dimension~$n-2$, where $ \dim W =
  2n$.
\end{propo}
\begin{proof}
  The standard way to treat bubbled curves consists in considering
  them as elements in a bubble tree: here such a tree is composed by a
  simple holomorphic disk $u_0\co (\Disk{2},\S^1) \to
  (W,\plastikstufe{S})$ and holomorphic spheres
  $u_1,\dotsc,u_{k^\prime}\co\S^2\to W$.  These holomorphic curves
  are connected to each other in a certain way.  We formalize this
  relation by saying that the holomorphic curves are vertices in a
  tree, ie\ in a connected graph without cycles.  We denote the
  edges of this graph by $\{u_i,u_j\}$, $0\le i<j\le k^\prime$.

  Now we assign to any edge two nodal points $z_{ij}$ and $z_{ji}$,
  the first one in the domain of the bubble $u_i$, the other one in
  the domain of $u_j$, and we require that
  $\evaluationmap_{z_{ij}}(u_i) = \evaluationmap_{z_{ji}}(u_j)$.  For
  technical reasons, we also require nodal points on each holomorphic
  curve to be pairwise distinct.  To include into the theory, trees
  with more than one bubble connected at the same point to a
  holomorphic curve, we add ``ghost bubbles''.  These are constant
  holomorphic spheres inserted at the point where several bubbles are
  joined to a single curve.  Now all the links at that point are
  opened and reattached at the ghost bubble.  Ghost bubbles are the
  only constant holomorphic spheres we allow in a bubble tree.

  The aim is to give a manifold structure to these bubble trees.
  Unfortunately this is in general not always possible, because
  already for a single sphere, one can only obtain regularity of the
  Cauchy--Riemann operator, if the sphere is simple.

  Instead, we note that the image of every bubble tree is equal to the
  image of a simple bubble tree, that means, to a tree, where every
  holomorphic sphere is simple and any two spheres have different
  image.  Since we are only interested in the image of the evaluation
  map on the bubble trees, it is for our purposes equivalent to
  consider the simple bubble tree instead of the original one.  The
  disk $u_0$ is always simple, and does not need to be replaced by
  another simple curve.

  Let $u_0,u_1,\dotsc,u_{k^\prime}$ be the holomorphic curves
  composing the original bubble tree, and let $A_i\in H_2(W)$ be the
  homology class represented by the holomorphic sphere $u_i$.  The
  simple tree is composed by $u_0,v_1,\dotsc,v_{k}$ such that for
  every $u_j$ there is a bubble sphere $v_{i_j}$ with
  \begin{align*}
    u_j(\S^2) &= v_{i_j}(\S^2)
  \end{align*}
  and in particular $\smash{A_j = m_jB_{i_j}}$, where $\smash{B_{i_j} = [v_{i_j}]\in
  H_2(W)}$ and $m_j\ge 1$ is an integer.  Denote the sum
  $\smash{\sum_{j=1}^{k^\prime}A_j}$ by $A$ and the sum $\smash{\sum_{i=1}^{k}B_i}$ by $B$.
  Below we will now compute the dimension of this simple bubble tree.

  The initial bubble tree $\smash{u_0, u_1,\dotsc, u_{k^\prime}}$ is the limit
  of a sequence in the moduli space
  $\mathcal{M}(W,\plastikstufe{S},z_0)$.  Hence the connected sum
  $u_\infty := u_0\sharp \dots \sharp u_{k^\prime}$ is, as element of
  $\pi_2(W,\plastikstufe{S})$, homotopic to a disk $u$ in the bishop
  family, and the Maslov indices
  \begin{equation*}
    \mu(u) :=
    \mu(u^*TW,u^*T\plastikstufe{S})
    \quad\text{and}\quad
    \mu(u_\infty) :=
    \mu(u_\infty^*TW,u_\infty^*T\plastikstufe{S})
  \end{equation*}
  have to be equal.  By \fullref{maslov index der bishop kurven}
and standard rules for the Maslov index, we obtain
  $$
    2 = \mu(u) = \mu(u_\infty) = \mu(u_0) + \sum_{j=1}^{k^\prime}
    2c_1([u_j]) = \mu(u_0) + 2c_1(A).
  $$
  The dimension of the unconnected set of holomorphic curves
  $$\mathcal{M}_{[u_0]}(W,\plastikstufe{S},z_0) \times \prod_{j=1}^k
  \mathcal{M}_{B_j}(W)$$ for the simple bubble tree is
  \begin{align*}
    \bigl(n + \mu(u_0)\bigr) + \sum_{j=1}^k 2\,\bigl(n+c_1(B_j)\bigr)
    & = n + 2 -
    2c_1(A) + 2nk + \sum_{j=1}^k 2c_1(B_j) \\
    & = n + 2 + 2nk + 2\, \bigl(c_1(B)-c_1(A)\bigr).
  \end{align*}
  In the next step, we want to consider the subset of connected
  bubbles, ie\ we choose a total of $k$ pairs of nodal points, which
  then have to be pairwise equal under the evaluation map.  The nodal
  points span a manifold
  \begin{align*}
    Z(2k) &\subset \bigl\{(1,\dotsc,2k) \to
    \Disk{2}\dot\cup\,\S^2\dot\cup \dots \dot\cup\,\S^2\bigr\}
  \end{align*}
  of dimension $4k$.  The dimension reduction comes from requiring
  that the evaluation map
  \begin{align*}
    \evaluationmap\co \mathcal{M}_{[u_0]}(W,\plastikstufe{S},z_0)
    \times \prod_{j=1}^k \mathcal{M}_{B_j}(W) \times Z(2k) \to W^{2k}
  \end{align*}
  sends pairs of nodal points to the same image in the symplectic
  manifold.  By regularity and transversality of the evaluation map to
  the diagonal submanifold $\triangle(k)\hookrightarrow W^{2k}$, the
  dimension of the space of holomorphic curves is reduced by the
  codimension of $\triangle(k)$, which is $2nk$.

  As a last step, we have to take the quotient by the automorphism
  group to obtain the moduli space.  The dimension of the automorphism
  group is $6k + 2$ (the last term corresponds to the automorphism
  group of the holomorphic disk with one marked point on its
  boundary).  Hence the dimension of the total moduli space is
  \begin{multline*}
    n + 2 + 2nk + 2\, \bigl(c_1(B)-c_1(A)\bigr) +
    4k - 2nk - (6k +2) \\
    = n - 2k + 2\, \bigl(c_1(B)-c_1(A)\bigr) \ge n - 2k.
  \end{multline*}
  The inequality holds because by the assumption of semipositivity,
  all the Chern classes are nonnegative on holomorphic spheres, and
  all coefficients $n_j$ in the difference $c_1(B)-c_1(A) = \sum_j
  c_1(B_j) - \sum_i c_1(A_i) = \sum_j c_1(B_j) - \sum_i m_i
  c_1(B_{j_i}) = \sum_j n_j c_1(B_j)$ are nonpositive integers.
\end{proof}

\section{Applications}

\subsection[Exotic contact structures on (2n-1)-dimensional Euclidean space]{Exotic contact structures on $\R^{2n-1}$}
\label{sec: exotische kontakstruktur}

A contact structure on $\R^{2n-1}$ containing a plastikstufe is
obviously exotic, because it cannot embed into the standard sphere
$(\S^{2n-1},\alpha_0)$.  This application is due to Chekanov
and Gromov \cite[2.4.$D_2^\prime$.(c)]{Gromov_Kurven}.

Instead of using the most general setup, we will just give one
example.  Bates and Peschke \cite{BatesPeschke} have
constructed an exotic symplectic structure $\omega = -d\lambda$ (see
also \hbox{\cite[Example~13.8]{McDuffSalamonIntro}}) on $\R^4$ that contains
a Lagrangian torus $\T^2$ such that $\left.\lambda\right|_{\T^2}
\equiv 0$.  Let $\alpha_-$ be an overtwisted contact structure on
$\R^3$.  Then
\begin{align*}
  \bigl(\R^7 = \R^3\times \R^4, \alpha_- + \lambda \bigr)
\end{align*}
is an exotic contact structure, because it contains the plastikstufe
$\plastikstufe{\T^2}$.  To the author's knowledge, there is until now,
no other way to distinguishing these contact structures from the
standard one: classical invariants fail to do so, and contact homology
on open manifolds has not yet been rigorously developed.

Other examples of exotic contact structures on $\R^{2n-1}$ were known
before (eg\ Muller shows that the symplectization of an
exotic $\R^5$ she constructs, contains a Lagrangian sphere
\cite{MullerExoticR6}), but with the plastikstufe it is easy to
construct examples in any dimension (embed the neighborhood of a
plastikstufe into $\R^{2n-1}$, and use the $h$--principle to extend
the contact structure over the whole Euclidean space).

\subsection{Filling with holomorphic curves in higher
  dimensions}\label{sec: filling}

In this section, we will apply the ideas in the proof of
\fullref{hauptsatz} to an example not directly related to the
plastikstufe.  This example can be seen as a generalization of
\cite[2.4.$\mathrm{D}_1^\prime$]{Gromov_Kurven} to higher dimensions,
and the main difficulty consists again in finding a replacement for
positivity of intersections.

Consider a closed contact $(2n-1)$--manifold $(M,\alpha)$, which is
the convex boundary of a semipositive symplectic manifold
$(W,\omega)$.  Any Darboux chart $U\subset M$ contains subsets 
contactomorphic to $\S^{2n-1}-\{(0,\dotsc,0,1)\}$ with the standard
contact form
\begin{align*}
  \alpha_0 &= \sum_{j=1}^n \bigl(x_j\,dy_j - y_j\,dx_j\bigr),
\end{align*}
where $z_j = x_j + iy_j$ are complex coordinates of $\C^n$, and
$\S^{2n-1}$ is embedded in the standard way into $\C^n$.  Note that
the canonical $\SO(n)$--action on $\C^n$ (ie\ the one by matrix
multiplication) restricts to the sphere and leaves $\alpha_0$
invariant.

Let $P_{t_0}$ be the $3$--plane
\begin{align*}
  P_{t_0} &= \bigl\{(x+iy,z+it_0,0,\dotsc,0)\bigm|
  x,y,z\in\R\bigr\}.
\end{align*}
Denote the intersection $P_{t_0} \cap \S^{2n-1}$ by
$\S^2_{t_0}$.  It is a $2$--sphere (if $t_0\in (-1,1)$), and the
foliation induced by $\alpha_0$ has only two singularities, at the
north pole $N = (0,\smash{\sqrt{\smash{1-t_0^2}\vphantom{^h}}} + it_0,0,\dotsc,0)$ and the south pole $S = (0,-\smash{\sqrt{\smash{1-t_0^2}\vphantom{^h}}} + it_0,0,\dotsc,0)$.
Every leaf circles down from $N$ to $S$.

Now we consider the flow-out $F_{t_0} := \SO(n-1)\cdot \S^2_{t_0}$
obtained by taking the set of all $\SO(n-1)$--orbits, where we use the
embedding
\begin{align*}
  \SO(n-1) \hookrightarrow \SO(n)\\
  A \mapsto
  \begin{pmatrix}
    1 & 0 \\ 0 & A
  \end{pmatrix}.
\end{align*}
Each orbit is diffeomorphic to an $(n-2)$--sphere $\S^{n-2} \cong
\SO(n-1)/\SO(n-2)$, and lies inside the contact structure
$\ker\alpha_0$.  Hence it follows that $F_{t_0} \cong \S^2\times
\S^{n-2}$ is a maximally foliated submanifold, with the two elliptic
singularities $\{S\}\times \S^{n-2}$ and $\{N\}\times \S^{n-2}$.  The
regular leaves connect the upper singularity with the lower one.  We
just have shown that we can find such a maximally foliated submanifold
$F_{t_0}$ in any Darboux chart.

By choosing now the almost complex structure from \fullref{sec:
  bishop family} around both elliptic singularities, and extending
this to a generic $\omega$--compatible structure on $W$, we obtain
around each of the poles of the sphere a Bishop family, and the aim
will be to show that in fact the upper and lower family lie in the
same moduli space.

In the same way, as in the proof of \fullref{hauptsatz}, it can be
excluded that the moduli space has any boundary components apart from
the Bishop ends, so the moduli space is either connected with two
boundary components, one for each of the Bishop families, or it
consists of two disconnected spaces, each with one boundary component
(see \hbox{\fullref{bild: moduli raum mit zwei enden}}).  If the moduli
space was disconnected, then the Bishop end would represent a trivial
homology class, but its image under the homology class is not, giving
a contradiction.

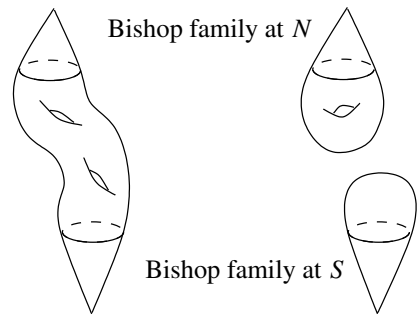
\begin{figure}[ht!]
  \begin{center}
    \begin{picture}(0,0)%
      \includegraphics{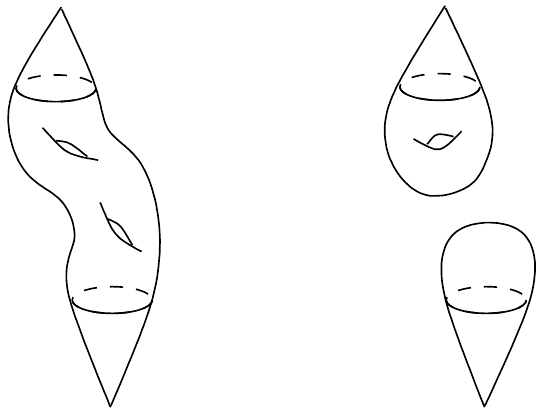}%
    \end{picture}%
    \setlength{\unitlength}{4144sp}%
    \begingroup\makeatletter\ifx\SetFigFont\undefined%
    \gdef\SetFigFont#1#2#3#4#5{%
      \reset@font\fontsize{#1}{#2pt}%
      \fontfamily{#3}\fontseries{#4}\fontshape{#5}%
      \selectfont}%
    \fi\endgroup%
    \begin{picture}(2480,1859)(2053,-7804)
      \put(2881,-7576){\makebox(0,0)[lb]{\smash{{\SetFigFont{8}{9.6}{\familydefault}{\mddefault}{\updefault}{\color[rgb]{0,0,0}Bishop family at $S$}%
            }}}}
      \put(2656,-6136){\makebox(0,0)[lb]{\smash{{\SetFigFont{8}{9.6}{\familydefault}{\mddefault}{\updefault}{\color[rgb]{0,0,0}Bishop family at $N$}%
            }}}}
    \end{picture}%
    \caption{The moduli space has either a single component with two
      ends or two components each with a single end.}\label{bild:
      moduli raum mit zwei enden}
  \end{center}
\end{figure}

\appendix

\section{Outlook and open problems}

The notion of overtwistedness plays a very central role in
$3$--dimensional contact topology.  The main implications are, as
mentioned in the introduction, the nonfillability and the easy
classification of such manifolds.  We have tried to generalize the
definition to higher dimensions, and have proved that our definition
implies nonfillability, but clearly, the easy classification of
overtwisted $3$--manifolds is for contact topology the more important
of the two properties.  In particular, it has been shown that any
oriented $3$--manifold supports overtwisted contact structures.  Thus
it is very disappointing that we have not been able to find a single
closed contact manifold of higher dimensions containing a
plastikstufe.

An interesting future research goal could consist in trying to find a
relation between the plastikstufe and the work of Giroux.  He
discovered that contact manifolds can be decomposed into open books,
which are compatible with the contact structure.  He also recognized
that any contact open book, whose monodromy is composed only by
right-handed Dehn twists, leads to a fillable contact manifold.
Guided by this realization, one could try to find for example a
plastikstufe in the sphere $(\S^{2n+1},\alpha_-)$, whose open book
decomposition has page $P \cong (T^*\S^n, d\lcan)$ and whose monodromy
consists of a single left-handed Dehn twist.  In the $3$--dimensional
case, it is easy to find explicitly an overtwisted disk.  More
generally, if a plastikstufe could be found, one could ask:

\begin{question}
  Can one read off from an open book decomposition, whether the
  contact manifold, contains a plastikstufe?
\end{question}

Finding examples of closed manifolds is the most immediate problem,
but other questions are also interesting.  The only application of the
plastikstufe so far is the detection of exotic contact structures on
$\R^{2n-1}$.  There are many constructions leading to such manifolds,
which could potentially lead to nonequivalent contact forms.  Use for
example different overtwisted contact structures on $\R^3$ in
\fullref{sec: exotische kontakstruktur} to create exotic
structures on~$\R^7$.

\begin{question}
  Can one somehow distinguish some of the exotic contact structures on
  $\R^{2n-1}$ containing a plastikstufe?
\end{question}

In dimension~$7$ and higher, there are many different plastikstufes,
because one can use different choices for the singular set~$S$.

\begin{question}
  Are the different versions of the plastikstufe in higher dimension
  somehow equivalent?
\end{question}

\section[The almost complex structure on the cotangent
  bundle]{The almost complex structure on the cotangent
  bundle}\setobjecttype{App}\label{sec: cotangent bundle}

The following statements about the cotangent bundle can certainly be
found in many references, but for completeness, we still repeat them
here: the aim will be to associate to any Riemannian manifold $(M,g)$
a natural metric on the cotangent bundle $T^*M$ (ie\ a bundle metric
on $T(T^*M)$), and use this to choose an almost complex structure $J_g$
on $T^*M$.

Let $\pi\co T^*M \to M$ be the standard projection, and let $g^\dagger$
be the bundle metric on $T^*M$ induced by $g$.  The vertical bundle
$V(T^*M)$ of $T(T^*M)$, ie\ the kernel of $\pi_*\co T(T^*M) \to TM$,
can be naturally identified with the bundle $\pi^*(T^*M)$ over $T^*M$
by taking two covectors $\beta_1,\beta_2\in \smash{T^*_p M}$, and considering
the derivative of the path $\beta_1 + t\,\beta_2\subset T^*_pM$ at
$t=0$.  This identification makes it possible to use $\smash{g^\dagger}$ to
define a bundle metric on $V(T^*M)$.

The Levi-Civita connection gives the natural splitting
\begin{align*}
  T(T^*M) &= H(T^*M) \oplus V(T^*M)
\end{align*}
into horizontal and vertical bundle.  Denote the vertical part of a
vector $v\in T(T^*M)$ by $v_V$ and the horizontal one by $v_H$.  This
splitting induces a natural metric $\smash{\tilde g}$ on $T^*M$ (ie\ a
bundle metric on $T(T^*M)$)
\begin{align*}
  \tilde g(v,w) &:= g(\pi_* v, \pi_* w) + g^\dagger(v_V,w_V),
\end{align*}
where we used the natural identification of the vertical bundle
described above.

\begin{theorem}\label{laenge auf faser J und lcan}
  There is a unique almost complex structure $J_g$ on $T^*M$ that is
  compatible with $d\lcan$ and $\tilde g$.  Furthermore the function
  \begin{align*}
    F(\POS,\MOM) &:= \frac{1}{2}\,g^\dagger(\MOM,\MOM)
  \end{align*}
  on $T^* M$ satisfies
  \begin{align*}
    dF \circ J_g &= -\lcan.
  \end{align*}
\end{theorem}
\begin{proof}
  We need to check that there is a unique solution $J_g$ for the
  equation
  \begin{align*}
    d\lcan(\cdot,J_g\cdot) &= \tilde g(\cdot,\cdot),
  \end{align*}
  such that $d\lcan(J_g\cdot,J_g\cdot) = d\lcan(\cdot,\cdot)$ and
  $J_g^2 = -\1$.

  The equations are independent of any chart, hence it suffices to
  check them at every point in one special chart explicitly.  Choose
  for a point $\POS_0\in M$ a geodesic normal chart, ie\ coordinates
  $(q^1,\dotsc,q^n)$ such that the $\langle
  \partial_{q^1},\dotsc,\partial_{q^n}\rangle$ form at $\POS_0$ an
  orthonormal basis, and such that all Christoffel symbols vanish at
  that point.  For the cotangent bundle, we obtain a chart
  $(q^1,\dotsc,q^n,p_1,\dotsc,p_n)$, where the vertical bundle
  \begin{align*}
    V_{(\POS_0,\MOM)}(T^*M) &= \Bigl\{(\POS_0,\MOM;0,\dot \MOM)\Bigm|
    \MOM,\dot\MOM \in \R^n\Bigr\},
  \end{align*}
  and horizontal bundle
  \begin{align*}
    H_{(\POS_0,\MOM)}(T^*M) &= \Bigl\{(\POS_0,\MOM;\dot \POS,0)\Bigm|
    \MOM,\dot\POS \in \R^n\Bigr\}
  \end{align*}
  at $\smash{T^*_{\POS_0}M}$ take the very easy form written above.  The
  metric $\tilde g$ can be represented at $\POS_0$ in the chart by the
  matrix $\smash{\left(\begin{smallmatrix} \1 & 0 \\0 & \1
    \end{smallmatrix}\right)}$, and the $2$--form $d\lcan$ by
  $\smash{\left(\begin{smallmatrix} 0 & \1 \\ -\1 & 0
    \end{smallmatrix}\right)}$.  It follows that at the given point,
  $J_g$ is the map that sends $\smash{\partial_{q^j}}$ to $\smash{\partial_{p_j}}$ and
  $\smash{\partial_{p_j}}$ to $\smash{-\partial_{q^j}}$.  This solves the first claim
  of the theorem.

  To check the equality for the function $F$, we will again use a
  normal geodesic chart around $\POS_0$ as explained above.  A short
  computation at $\POS_0$ shows that (the $\smash{g^{ij}}$ denote the
  coefficients of the metric $g^\dagger$)
  \begin{align*}
    dF &= \frac{1}{2}\frac{\partial g^{ij}}{\partial
      q^k}\,p_ip_j\,dq^k + g^{ij} p_i\,dp_j = p_j\,dp_j.
  \end{align*}
  Since $J_g$ sends $\partial_{q^j}$ to $\partial_{p_j}$, it follows that
  \begin{align*}
    dF\circ J_g &= p_j\,dp_j\circ J_g = p_j\,dq^j = - \lcan
  \end{align*}
  at the considered point, but since there is a geodesic normal chart
  around any point, the equation holds everywhere.
\end{proof}

\bibliographystyle{gtart}
\bibliography{link}

\end{document}